\documentclass{article}
\usepackage{arxiv}
\usepackage{amsfonts}
\usepackage{amssymb}
\usepackage{amsthm}
\usepackage{bbm}
\usepackage{bm}
\usepackage{color}
\usepackage{hyperref}
\usepackage{graphicx}
\usepackage{physics}
\usepackage{systeme}
\usepackage{svg}
\usepackage{textgreek}
\usepackage{tikz}
\usetikzlibrary{positioning, fit}
\usepackage{amsfonts} 
\usepackage[normalem]{ulem}
\usepackage[utf8]{inputenc}
\newcommand{\vect}[1]{\boldsymbol{#1}}
\usepackage{booktabs} 
\usepackage{amsmath}

\title{Learning Effective SDEs \\ from Brownian Dynamics Simulations\\ of Colloidal Particles}

\author{
  Nikolaos Evangelou\\
 Johns Hopkins University\\
  Baltimore, MD 21218, USA\\
  \And
  Felix Dietrich\\
  Technical University of Munich\\
  Munich, 80333, Germany\\
  \And
  Juan M. Bello-Rivas\\
 Johns Hopkins University\\
  Baltimore, MD 21218, USA\\
  \And
  Alex Yeh\\
 Johns Hopkins University\\
  Baltimore, MD 21218, USA\\
  \And
  Rachel Stein \\
  Johns Hopkins University\\
  Baltimore, MD 21218, USA\\
  \And
  Michael A. Bevan \\
 Johns Hopkins University\\
  Baltimore, MD 21218, USA\\
  \And
  Ioannis G. Kevrekidis*\\
 Johns Hopkins University\\
  Baltimore, MD 21218, USA\\
  yannisk@jhu.edu
}

\begin{document}
\maketitle
\begin{abstract}
{We construct a reduced, data-driven, parameter dependent effective Stochastic Differential Equation (eSDE) for electric-field mediated colloidal crystallization using data obtained from Brownian Dynamics Simulations. We use Diffusion Maps (a manifold learning algorithm) to identify a set of useful latent observables. In this latent space we identify an eSDE  using a deep learning architecture inspired by numerical stochastic integrators and compare it with the traditional Kramers-Moyal expansion estimation. We show that the obtained variables and the learned dynamics  accurately encode  the  physics of  the  Brownian Dynamic Simulations. We further illustrate that our reduced model captures the dynamics of corresponding experimental data. Our dimension reduction/reduced model identification approach can be easily ported to a broad class of particle systems dynamics experiments/models.}
\end{abstract}

\section{Introduction}

The identification of nonlinear dynamical systems from experimental time series and image series data became an important research theme in the early 1990s \cite{krischer1993model,rico1992discrete,rico1993discrete}. After lapsing for almost two decades, it is now experiencing a spectacular rebirth.
A key element of the older work was the use of neural architectures \cite{GONZALEZGARCIA1998S965,rico1992discrete} (recurrent, convolutional, ResNet) motivated by traditional numerical analysis algorithms.
Importantly, such architectures allow researchers to identify {\em effective}, coarse-grained, mean-field type evolution models from {\em fine-scale} (atomistic, molecular, agent-based) data \cite{liu2015equation,imapd}.

\begin{figure}[ht]
\begin{center}
\includegraphics[scale=0.25]{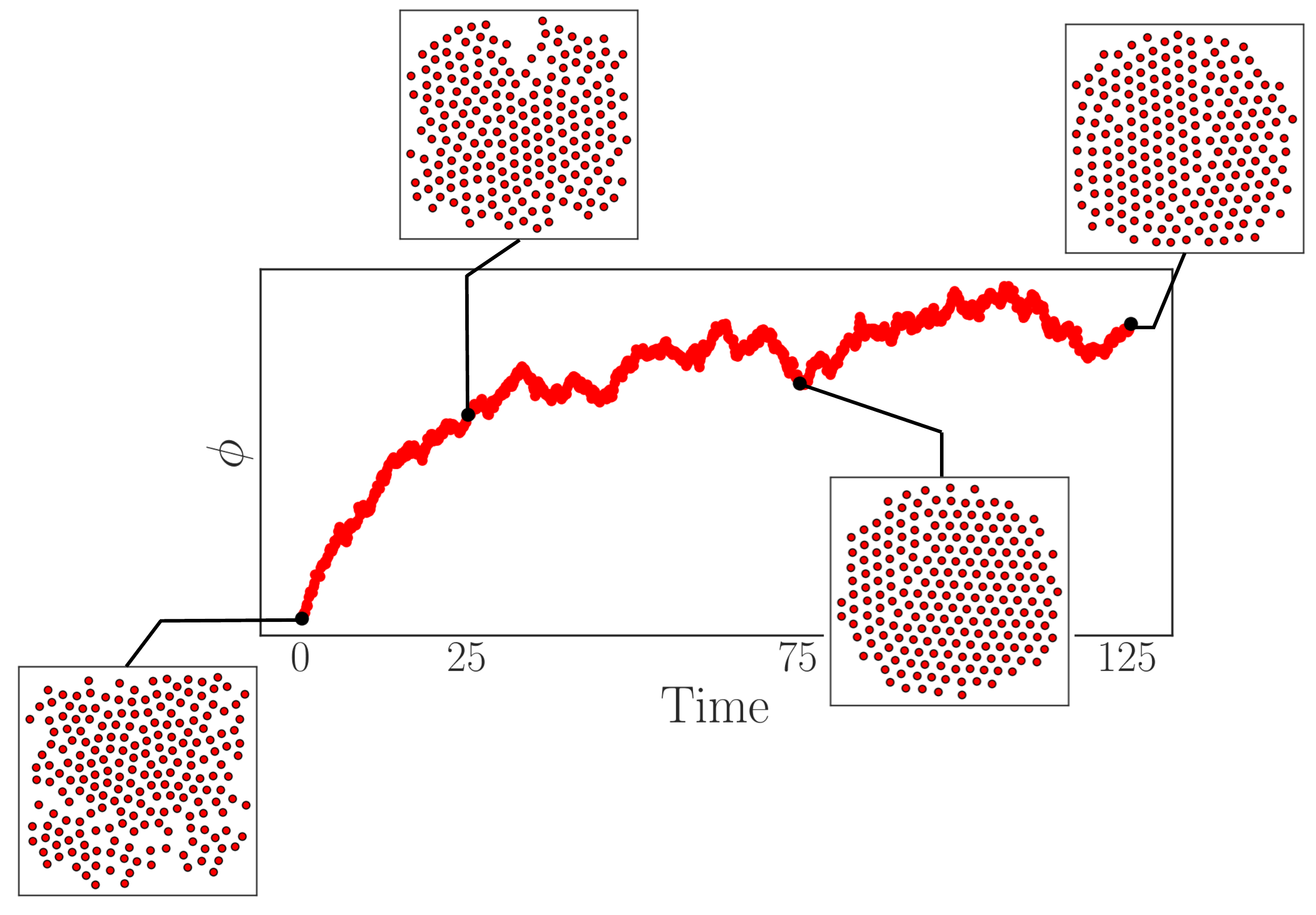}
\caption{A trajectory of an effective, reduced eSDE in the data driven collective coordinate $\phi$ for electric-field mediated colloidal crystallization.} 
\label{fig:Caricature_SDE}
\end{center}
\end{figure}

In this paper, we identify coarse-grained, effective stochastic differential equations (eSDE) for colloidal particle self-assembly based onfine-grained, Brownian dynamics simulations under the influence of electric fields \cite{yang2016dynamic,edwards2013size}. We demonstrate that the identified eSDE encodes accurately the physics of the Brownian Dynamic simulations and captures the dynamics of corresponding experimental data. Those experiments have previously been shown to quantitatively match to BD simulations {\em at equilibrium} in terms of time-averaged distribution functions \cite{edwards2013size,juarez2009interactions,juarez2011kt}.
Figure~\ref{fig:Caricature_SDE} shows a sample path of a latent space trajectory $\left\lbrace t,\phi(t)\right\rbrace_{t\geq 0}$ computed through our learned eSDE. The corresponding instantaneous particle conformations are indicated at representative points along the trajectory.
A key feature of our work is the selection of the coarse-grained observables (the variables of our eSDE) in a data-driven manner, using manifold learning techniques like Diffusion Maps \cite{COIFMAN20065}. The dependence of the dynamics on physical control parameters (here a driving voltage) is included in the neural architecture and learned during training. 
A second key feature is that the neural network architecture for eSDE identification is not based on established Kramers-Moyal estimation techniques e.g. \cite{gradivsek2000analysis,risken1996fokker},
but rather (in the spirit of the early work mentioned above) on numerical stochastic integration algorithms \cite{dietrich2021learningsdes}. 

The motivation of the effective SDE is to better understand the ability to assemble nano- and micro- colloidal particles into ordered materials and controllable devices. This could provide the basis for emerging technologies (e.g. photonic crystals, meta-materials, cloaking devices, solar cells, etc.), \cite{arpin2010multidimensional} but also impact traditional applications (e.g. ceramics, coatings, minerals, foods, drugs \cite{russel1989formulation},\cite{zukoski1995particles}). Despite the range of applications employing microscopic colloidal particles, current state-of-the-art \cite{hendley2021anisotropic} capabilities for manipulating microstructures in such systems are limited in two ways: (a) the degree of order that can be obtained, and (b) the time required to generate ordered structures. Both of these limitations are due to fundamental problems with designing, controlling, and optimizing (i.e. engineering) the thermodynamics and kinetics of colloidal assembly processes. 
By learning effective collective models, our approach in this work paves the way towards data-driven model-based control of the kinetics of  colloidal assembly processes using parametric driving (e.g. via an electric field) \cite{juarez2012feedback}. Crucially, we also tackle the issue of interpretability of the learned effective dynamic model by exploring relations between data-driven and candidate physically meaningful observables. Those coarse physical observables are \textit{order parameters} that provide intuition for the colloidal self-assembly process \cite{tang2016optimal,tang2017construction,zhang2020controlling}.
We combine the data-driven detection of effective latent spaces with the neural network based, numerical analysis inspired, identification of parameter-dependent stochastic eSDEs with state-dependent diffusion. This is based on fine scale data from both Brownian dynamics simulations and from experimental colloidal crystallization movies, and the results are compared.

Developing low-dimensional surrogate models for physical systems has been explored by a number of authors. We report some approaches that utilize machine learning and/or dimensionality reduction here that could be beneficial to the reader. The authors in \cite{kopelevich2005coarse} identified an effective, coarse grained Fokker-Planck using Kramers-Moyal with an application to micelle-formation of surfactant molecules. The identified equation in \cite{kopelevich2005coarse} was constructed in terms of the physical coarse variable, size of the cluster of the surfactant molecules. The authors in \cite{beltran2011smoluchowski} constructed an 1D Smoluchowski equation in terms of coarse physical variables (radius of gyration or the average crystallinity) for small colloidal systems of 32-particles.
The authors in \cite{coughlan2019non}  used simulation and experimental colloidal ensembles with smaller than 14 particles to fit two dimensional Fokker-Planck and Langevin equations. The two coarse variables in which the dynamics are being identified capture the condensation and anisotropy of those small ensembles. 
A detailed review that summarizes applications of machine learning to discover collective variables and for sampling enhancement was conducted by the authors \cite{sidky2020machine}.
A framework to advance the simulation time by learning the effective dynamics (LED) of molecular systems was proposed by the authors in \cite{vlachas2021accelerated}. 
LED uses mixture density network (MDN) autoencoders to learn a mapping between the molecular systems and latent variables and evolves the dynamics using long short-term memory MDNs.
In the context of accelerating molecular simulations, the authors in \cite{fu2022simulate} proposed a framework tested on polymeric systems that utilize graph clustering to obtain coarse observable and allows to model system's evolution for long-time dynamics.

The main machine learning tools of our work involve (a) utilizing a dimensionality reduction scheme that discovers a lower dimensional structure of a given data set and (b) a deep neural network architecture that learns an eSDE.

Regarding the first aspect of dimensionality reduction a wide range of techniques have been proposed for discovering a set of reduced observables. Among others, Principal Component Analysis \cite{pearson1901liii}, Isomap \cite{tenenbaum2000global}, Local Linear Embedding \cite{roweis2000nonlinear}, Laplacian Eigenmaps \cite{belkin2003laplacian}, Autoencoders \cite{kramer1991nonlinear} and our method of choice:  Diffusion Maps~\cite{COIFMAN20065}. The Diffusion Maps algorithms enables the discovery of reduced coordinates when data are sampled from signal processing \cite{talmon2013diffusion}, from networks \cite{rajendran2016data},   from (stochastic) differential equations \cite{NADLER2006113,Chiavazzo_Reduced_Models} but also from Molecular \cite{imapd} and Brownian simulations \cite{yang2016dynamic}. 

A traditional approach for learning eSDEs has been the Kramers-Moyal expansion \cite{gradivsek2000analysis,risken1996fokker,liu2015equation} and a detailed description of this approach is given in Section \ref{sec:Kramers_Moyal}. For non-Gaussian stochastic differential equations, modifications of the tradiational Kramers-Moyal expansion have also been proposed \cite{lu2021extracting,li2021data}. In \cite{oleary2021stochastic} the authors proposed a stochastic  physics-informed neural network framework (SPINN) that minimizes the distance between the predicted moments of the network (drift and diffusivity) from moments computed with Kramers-Moyal. The authors in \cite{boninsegna2018sparse} proposed an extension of the framework called Sparse Identification of Nonlinear Dynamics (SINDy) that can be used for stochastic dynamical systems. The authors in \cite{yang2020generative} proposed a physics informed generative model termed \textit{generative ensemble-regression}  that learns to generate \textit{fake} sample paths  from given densities at several points in time, without point-wise paths correspondence. The authors in \cite{li-2020} extracted an eSDE from long time series data in a memory-efficient way, including learning the eSDE in latent variables. This approach is valuable if the data is available as a few, long time series. In our approach we handle pairs of successive snapshots instead. The most similar approach to learning eSDEs to the one selected for our work is \cite{hasan2021identifying}. The authors introduce a Variational  Autencoder (VAE) framework for recovering latent dynamics governed by an eSDE. In their method, the latent space and the stochastic differential equation are identified together within the VAE scheme. Their loss function is also based on the Euler-Maruyama scheme.

Our work deviates from the approaches mentioned above in three key aspects: (a) we explicitly separate the latent space construction from learning the eSDE; (b) we extend the loss function informed by numerical integration schemes from \cite{dietrich2021learningsdes} to allow for additional parameter dependence.
Our latent space is defined through  Laplace-Beltrami operator eigenfunctions, so, different from \cite{hasan2021identifying}, (c) our latent space coordinates are invariant to isometry and sampling density in the original space by construction.

\section{Methodology}
\subsection{Brownian Dynamics}

\begin{figure}[ht]
    \centering
    \includegraphics[scale=0.28]{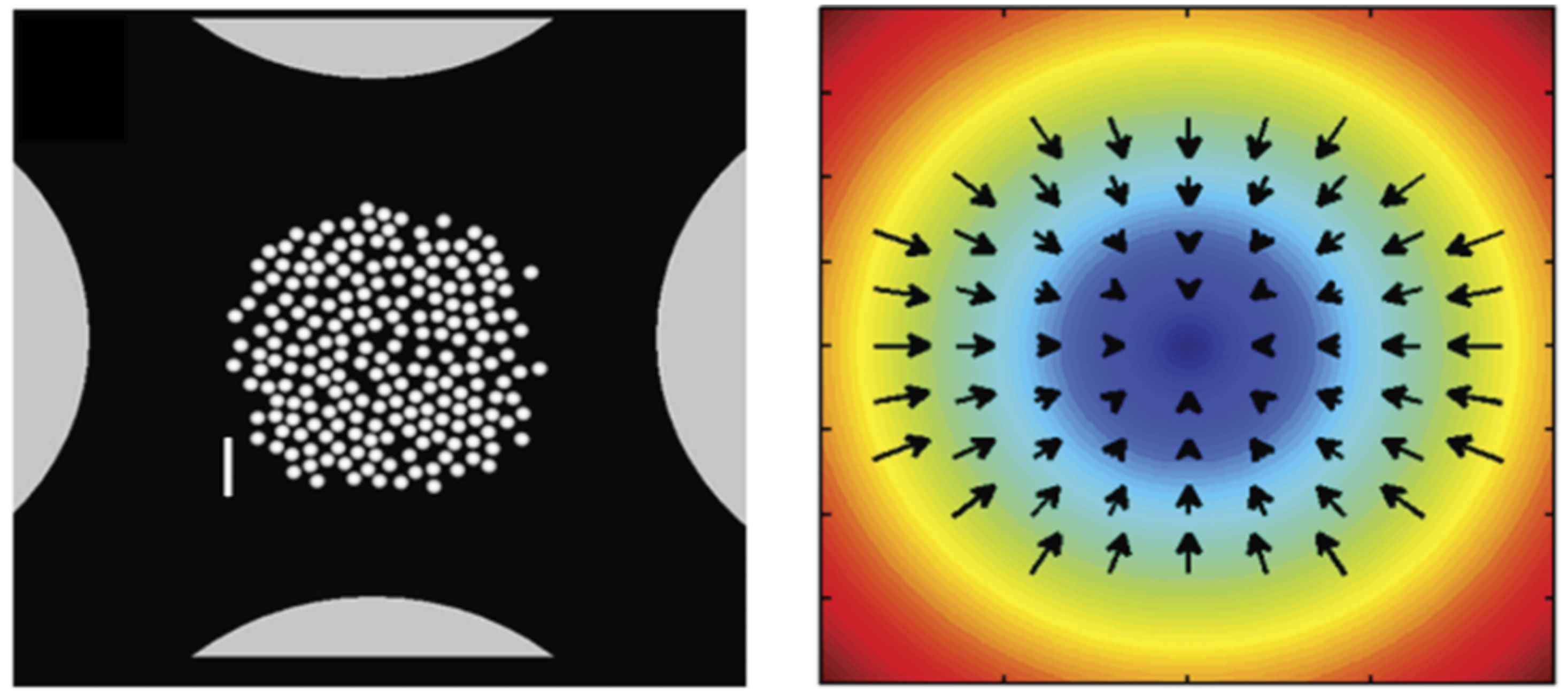}
    \caption{[Left]  Top view of simulated experiments of quasi 2D configurations of $N=210$ colloidal particles compressed with a quadrupole electrode. [Right] Electric field magnitude contour plot in the vicinity of the quadrupole electrode center. The arrows indicate the relative magnitude and direction of force due to dipole-field interactions. Taken from J.Chem. Phys \textbf{144}, 204904 (2016) with permission.}
    \label{fig:dynamic_colloidal}
\end{figure}

We model electric field-mediated quasi-2D colloidal assembly in the presence of a quadrupole electrode. An illustration of the set up is shown in Figure \ref{fig:dynamic_colloidal} In our simulations, each configuration consists of $N=210$ particles. The interactions between the colloidal particles are electrostatic double layer repulsion $u^{pp}_{e,i,j}$, dipole-field potentials $u^{pf}_{de,i}$ and dipole-dipole interaction potential $u_{dd,i,j}^{pp}$. The electrostatic repulsion,  $u^{pp}_{e,i,j}$, between two particles $i$ and $j$ is computed by.

\begin{equation}
\label{eq:electrostatic_repulsion}
    u_{e,i,j}^{pp}(r_{i,j}) = B^{pp} \exp {\left[-\kappa(r_{ij} - 2\alpha) \right]}.
\end{equation}

In Equation \eqref{eq:electrostatic_repulsion} $r_{ij}$ denotes the center-to-center distance between the patricles, $\alpha$ is the radius of each particle and $B^{PP}$ is the electrostatic repulsion pre-factor between colloidal particles.

The dipole field potential $u^{pf}_{de,i}$ in the spatially varying electric field for each particle $i$ is computed by
\begin{equation}
    u^{pf}_{de,i} (\mathbf{r}_i) = -2kT\lambda f_{cm}^{-1}[E(\mathbf{r}_i)/E_0]^2,
\end{equation}

where $\mathbf{r}_i$ is the position of the $i^{\text{th}}$ particle, $k$ is the Boltzmann's constant, $T$ is the temperature, $f_{cm}$ is the Clausius-Mossotti factor, $\lambda$ is a non-dimensional amplitude given by the relation $\lambda = \frac{\pi\epsilon_m \alpha^3(f_{cm}E_0)^2}{kT}$, $\epsilon_m$ is the medium dielectric constant, the local electric field magnitude is given by $E(\mathbf{r}_i$). The constant $E_0$ is given by
the expression
\begin{equation}
    E_0  = \frac{1}{\sqrt{8}}(V_{pp}/d_g)
\end{equation}

where $V_{pp}$ denotes the peak-to-peak voltage and $d_g$ the electrode gap. The dipole-dipole interaction potential $u_{dd,i,j}^{pp}$ between two particles $i$ and $j$ is estimated by
\begin{equation}
    u_{dd,i,j}^{pp}(\mathbf{r}_{ij}) = -kT\lambda P_2(cos\theta_{ij})(2\alpha/r_{ij})^3[E(\mathbf{r}_i/E_0)]^2.
    \label{eq:dipole-dipole-potential}
\end{equation}

$P_2(cos{\theta_{ij}})$ is the second Legendre polynomial, $\theta_{ij}$ denotes the angle between the particle centers and the electric field direction. 

The electric field at the center of the quadrupole can be approximated by the expression
\begin{equation}
    \Bigg|\frac{E(\mathbf{r}_i)}{E_0} \Bigg| = \frac{4r}{d_g}
\end{equation}

where $r$ is the distance from the quadrupole center.

The motion of the Brownian particles is governed by the equation
\begin{equation}
    \mathbf{r}(t + \Delta t) = \mathbf{r}(t) + \frac{\mathbf{D}^P}{kT}(\mathbf{F}^P + \mathbf{F}^B)\Delta t + \nabla \cdot \mathbf{D}^P \Delta 
 t
\end{equation}

where $ \langle\mathbf{F}^B \rangle =0$, $\langle \mathbf{F}^B(t_1)(\mathbf{F}^B(t_2)^\text{T} = 2(kT)^2(\mathbf{D}^P)^{-1}\delta(t_1 - t_2)$, $\mathbf{r}(t)$ denotes the position vector for all the $N$ particles 
at time $t$, $\mathbf{F}^B$ denotes the Brownian force vector and $\mathbf{F}^P$ the total conservative force vector. The conservative force acting on each particle $i$ is given by
\begin{equation}
    \mathbf{F}_i^P = \nabla_{\mathbf{r}_i}\Big[u_{de,i}^{pf} + \sum_{j \neq i} (u_{e,i,j}^{pp} + u_{dd,i,j}^{pp})  \Big].
\end{equation}
$\mathbf{D}^P$ denotes the diffusivity tensor estimated by the Stoke-Einstein relation
\begin{equation}
    \mathbf{D}^P = kT(\mathbf{R}^P)^{-1}
\end{equation}
where $\mathbf{R}^P$ is the grand resistance tensor $\mathbf{R}^P$ given by
\begin{equation}
    \mathbf{R}^P = (\mathbf{M}^{\infty})^{-1} + \mathbf{R}_{2B} - \mathbf{R}^{\infty}_{2B}
\end{equation}

where $\mathbf{R}_{2B}$ are the pairwise lubrication interactions  and $(\mathbf{M}^{\infty})^{-1} -
\mathbf{R}^{\infty}_{2B}$ the many-bodied far-field interaction above a no-slip plane. All the parameters used for the BD simulations are included in Table 1 of the SI.

\subsection{Diffusion Maps}
\label{sec:Diffusion_Maps}
Introduced by \cite{COIFMAN20065}, Diffusion Maps offer a parametrization of a data set of points $\mathbf{X}=\{\vect{x}_i \}_i^N$ sampled from a manifold $\mathcal{M}$, where $\vect{x}_i \in \mathbb{R}^m$ by uncovering its \textit{intrinsic geometry}. This parametrization can then be used to achieve dimensionality reduction of the data set. This is obtained by initially constructing an affinity matrix $\mathbf{A}\in \mathbb{R}^{N \times N}$ through a kernel function, for example the Gaussian Kernel
\begin{equation}
    A_{ij} = \text{exp}\bigg(\frac{-||\vect{x}_i - \vect{x}_j ||^2}{2\varepsilon} \bigg),
\end{equation}
where $||\cdot||$ denotes a norm of choice. In this work we choose the $l^2$ norm; $\epsilon$ is a hyperparameter regulating the rate of decay of the kernel. To achieve a parametrization of $\mathbf{X}$ regardless of the sampling density, a normalization of $\mathbf{A}$ is performed as follows
\begin{equation}
    P_{ii} = \sum_{j=1}^N{A_{ij}},
\end{equation}
\begin{equation}
    \Tilde{\mathbf{A}} = \mathbf{P}^{-\alpha}\mathbf{A}\mathbf{P}^{-\alpha}
\end{equation}
where $\alpha =1$ is set to factor out the effect of sampling density.
The kernel $\mathbf{\Tilde{\mathbf{A}}}$ is further normalized
\begin{equation}
    W(\vect{x}_i,\vect{x}_j) = \frac{ \Tilde{A}(\vect{x}_i,\vect{x}_j)}{\sum_{j=1}^N \Tilde{A}(\vect{x}_i,\vect{x}_j)}
\end{equation}
so that the matrix $\mathbf{W}$ becomes a row stochastic matrix. 
The eigendecomposition of $\mathbf{W}$ results in a set of eigenvectors $\vect{\phi}$ and eigenvalues $\lambda$
\begin{equation}
    \textbf{W}\vect{\phi}_i = \lambda_i\vect{\phi}_i.
\end{equation}
To check if dimensionality reduction of \textbf{X} is possible, selection of the eigenvectors $\vect{\phi}$ that parameterize independent directions (non-harmonic eignvectors) is needed. In our work this selection was made by implementing the algorithm presented in \cite{dsilva2015parsimonious}. If the number of the non-harmonic eigenvectors is smaller than the original dimensions of $\mathbf{X}$ then Diffusion Maps achieves dimensionality reduction.

Our data ``points'', $\vect{x}_i$, consist of planar configurations of 210 particle locations obtained either from evolving computations or from experimental movies; our data set is $\mathbf{X} = \{\vect{x}_i \}_{i=1}^N$ where $\vect{x}_{i} \in \mathbb{R}^{210 \times 2}$. A number of preprocessing steps are performed before Diffusion Maps can be computed. All configurations are centered and aligned to a reference configuration by using Procrustes analysis, in particular the Kabsch algorithm \cite{kabsch1976solution,kabsch1978discussion}. The reference configuration was selected as the configuration that has the smallest value of the order parameter $Rg$ (see Section \ref{sec:order_parameters_energy}). Centering the data and applying the Kabsch algorithm removes the translational and rotational degrees of freedom. We then compute the density function, $\mathbf{f}_i$, for each configuration at the nodes of a grid and we normalize its integral to one. The density was estimated by a kernel density estimation using Gaussian Kernels in Python. More precisely, the \textit{gaussian\_kde} module from \textit{scipy} was used for this computation. The bandwidth for the kernel estimation was selected based on Scott's Rule \cite{scott2015multivariate}. The Diffusion Maps algorithm then is applied to the data set $\mathbf{F} =\{\mathbf{f}_i \}_i^N$ of the collected normalized density function discretizations. The density formulation $f_i$ eliminates the problem of permutational invariance of the particles in defining pairwise distances. As we mentioned also earlier the selection of the leading non-harmonic Diffusion Maps coordinates was made by the local linear algorithm proposed by the authors in \cite{dsilva2015parsimonious}. For our Diffusion Maps computations the \textit{datafold} package was used \cite{Lehmberg2020}.
\subsection{Nystr\"om Extension} 
Given a new out-of-sample data point, $x_\text{new} \notin \mathbf{X}$ (and subsequently $f_{new} \notin \mathbf{F})$, in order to embed it in the Diffusion Maps coordinates one might add it to the data set and recompute Diffusion Maps. However, this is computationally inefficient and will lead to a new Diffusion Maps coordinate system for every new point added in the data set. To avoid these issues the Nystr\"om Extension formula \cite{Nystrom_Extension,williams01_using_nystr} can be used
\begin{equation}
\phi_i(\vect{f_\text{new}}) = \frac{1}{\lambda_i}\sum_{j=1}^N \tilde{W}(f_{new},f_j)\phi_i(f_j),
\end{equation}
where $\phi_i(\vect{f_\text{new}})$ is the estimated value of the $i^{th}$ eigenvector for the new point $\vect{f_\text{new}}$, $\lambda_i$ is the corresponding eigenvalue, and $\phi_i(f_j)$ is the $j^{th}$ component of the $i^{th}$ eigenvector.

This formula is extremely useful in mapping trajectories either from the Brownian Dynamics simulations or from experimental snapshots to the Diffusion Maps coordinates (an operation called ``restriction''). Restricted long trajectories are used as a \textit{test} set to validate our estimated eSDEs.

\subsection{Learning SDEs from data}
In this section we describe two approaches to  \textit{estimate} SDEs from data.
Let $\vect{x}(t)$ be a stochastic vector-valued variable whose evolution is governed by the SDE
\begin{equation}
\label{eq:Langevin}
    d\vect{x}(t) = \vect{\nu}(\vect{x}(t))dt + \mathbf{\sigma}(\vect{x}(t))dB_t,
\end{equation}
where $\vect{\nu}:\mathbb{R}^m \to \mathbb{R}^m$ is the drift, $\mathbf{\sigma}:\mathbb{R}^m \to \mathbb{R}^{m \times m}$ is the diffusivity matrix, and $B$ a collection of $m$ one-dimensional Wiener processes. 
The dynamics of such process can be approximated by \textit{estimating} the two functions $\vect{\nu}$ and $\mathbf{\sigma}$. We show how this estimation can be performed, either from the statistical definition of the terms, based on the  Kramers-Moyal expansion \cite{risken1996fokker,gradivsek2000analysis}, or via a deep learning architecture inspired by stochastic numerical integrators \cite{dietrich2021learningsdes}.

\subsubsection{Kramers-Moyal expansion}
\label{sec:Kramers_Moyal}
For a stochastic process $\vect{x}(t)$, the differential change in time of its probability density $P(\vect{x},t)$  is given by
\begin{equation}
\label{eq:kramers_moyal}
    \frac{\partial P(\vect{x},t)}{\partial t}  = \sum_{n=1}^{\infty} \bigg(- \frac{\partial}{\partial \vect{x}}\bigg)^n \textbf{D}^{(n)}(\vect{x},t)P(\vect{x},t),
\end{equation}
which is known as the Kramers-Moyal expansion \cite{risken1996fokker}.
The moments of a transition probability, jumping from a position $x(t_k)$ to a nearby position $x(t_{k+h})$ in the next time step, are given by
\begin{equation}
\label{eq:kramers-moyal-coeff}
    D^{(n)}(\vect{x},t) = \frac{1}{n!}  \lim_{h \to 0} \frac{\langle[\vect{x}(t_{k+h} ) - \vect{x}(t_k)]^n \rangle}{h}. 
\end{equation}
where $\langle \cdot \rangle$ denotes the average. When $\vect{x}(t)$ is a Gauss-Markov process; only the first two moments of Equation \ref{eq:kramers_moyal} are non-zero,  and the Kramers-Moyal expansion reduces to the forward Fokker-Planck equation.
The Fokker-Planck equation provides an alternative description of the dynamics expressed by Equation \ref{eq:Langevin}. For $N$ variables, the Fokker-Planck equation is given by
\begin{equation}
    \frac{\partial \vect{P}}{\partial t} = - \sum_{i=1}^N \frac{\partial}{\partial x_i}\bigg(D_i^{(1)}(\vect{x})\vect{P} \bigg) + \sum_{i,j =1}^N \frac{\partial^2}{\partial x_i \partial x_j} \bigg( D_{ij}^{(2)}(\vect{x}) \vect{P} \bigg),
\end{equation}
where $\mathbf{D}^{(1)}$ and $\mathbf{D}^{(2)}$ are also the drift and diffusion coefficients and the connection with the coefficients of Equation $\eqref{eq:Langevin}$ is given by the expressions $
      \vect{\nu}(x) = \mathbf{D}^{(1)}, 
      \vect{\sigma}^2 = 2\mathbf{D}^{(2)}.
      $
The estimation of the drift and the diffusivity at a point $\vect{x}^i$ can be performed by multiple local parallel simulations (``bursts").

\begin{equation}
    \begin{aligned}
    &{\nu}_i(\vect{x}(t_k)) \approx \frac{1}{h}\langle {x}_i(t_{k+h})-{x}_i(t_k)\rangle, \\
    &{\sigma}^2_{ij}(\vect{x}(t_k)) \approx \frac{1}{h}\langle({x}_i(t_{k+h}) - \vect{x}_i(t_k))({x}_j(t_{k+h}) - \vect{x}_j(t_k))\rangle.
\end{aligned}
\end{equation}

\subsubsection{Deep Learning - Numerical Integrators}
\label{sec:Deep_Learning_sde}
The deep learning approach that we have followed for the identification of eSDEs is based on the work of \cite{dietrich2021learningsdes}. 
In this approach, the drift and diffusivity are estimated through two networks $\nu_\theta$ and $\sigma_\theta$, where $\theta$ are the weights of the networks. In our work we also introduce a small but meaningful modification of their method  
by also including a ``parameter neuron" along with the snapshot of inputs 
$\widetilde{\textbf{D}} = \{\vect{x}^{i}(t_{k+h}),\vect{x}^i(t_{k}), h^{i},p^{i} \}_{i=1}^N  $.
%
This modification allowed us to learn {\em parameter dependent eSDEs}. The parameter $p$ in our case is the applied  voltage to the particles ($V^*$).
The collected data required for this approach do not necessarily need to be sampled from long trajectories. Snapshots
$\widetilde{\textbf{D}}$
are sufficient as long as the region of interest is sampled densely enough. Here we introduce the scheme for the two-dimensional
case since our identified eSDE is also two-dimensional. Each snapshot $\widetilde{\vect{D}}^i$ in this network includes (a) a point at time $k$ in space  $\vect{x}^i(t_k) = (x_1^i(t_k),x_2^i(t_k))$; (b) its coordinates after a short time evolution  $\vect{x}^i(t_{k+h}) = (x_1^i(t_{k+h}),x_2^i(t_{k+h}))$; (c) the time interval between the two points, $h^i$; and (d) a parameter $p^i$ for the parameter dependent eSDE. Between different sampled snapshots the time step $h$ does not need to be uniform; in our case this property will prove to be quite useful as discussed in the results. 

The loss function used in our case (based on \cite{dietrich2021learningsdes}) is derived from the Euler-Maruyama scheme, a numerical integration method for SDEs. 
The scheme for the two-dimensional case,
\begin{multline}
\left\lceil
\begin{matrix}
 x_1^i(t_{k+h}) \\
  x_2^i(t_{k+h})
\end{matrix}
\right\rceil
=
\left\lceil
\begin{matrix}
 x_1^i(t_{k}) \\
  x_2^i(t_{k})
\end{matrix}
\right\rceil
+
h^i
\left\lceil
\begin{matrix}
 \nu_{\theta}(x_1^i(t_k),p^i) \\
  \nu_{\theta}(x_2^i(t_k),p^i)
\end{matrix}
\right\rceil
+ 
\left\lceil
\begin{matrix}
 \sigma_{\theta}(x_1^i(t_k),x_1^i(t_k),p^i) & \sigma_{\theta}(x_1^i(t_k),x_2^i(t_k),p^i) \\
  \sigma_{\theta}(x_2^i(t_k),x_1^i(t_k),p^i) &
  \sigma_{\theta}(x_2^i(t_k),x_2^i(t_k),p^i)
\end{matrix}
\right\rceil
\left\lceil
\begin{matrix}
 dB_{t_1} \\
  dB_{t_2}
\end{matrix}
\right\rceil
\end{multline}


where $dB_{t_1},dB_{t_2}$ are normally distributed around zero with variance $h^i$. This scheme has a similar form for higher dimensions.
This scheme implies that  each $\vect{x}^i(t_{k+h})$ is normally distributed,
\begin{equation}
\label{eq:distribution}
   \vect{x}^i(t_{k+h}) \sim \mathcal{N}\left(\vect{x}^i(t_k) + h^i\vect{\nu}_{\theta}(\vect{x}^i(t_k),p^i),h^i\mathbf{\sigma}_{\theta}(\vect{x}^i(t_k),p^i)^2\right)
\end{equation}
where the mean  $\vect{\mu}_{\theta}^i= \vect{x}^i(t_k) + h^i\vect{\nu}_{\theta}(\vect{x}^i(t_k),p^i)$ and the covariance matrix $\vect{\Sigma}^i_{\theta} = h^i\mathbf{\sigma}_{\theta}(\vect{x}^i(t_k),p^i)^2  $.

Under this assumption, we formalize a loss function that will lead to a maximization of the probability of Equation \eqref{eq:distribution}. This is achieved by combining the logarithm of the probability density of the multivariate normal distribution with the assumed mean and variance from Equation  \eqref{eq:distribution}:

\begin{multline}
    \mathcal{L}(\theta|\vect{x}^i(t_{k+h}),\vect{x}^i(t_k),h^i,p^i) :=   
   \log\big|\text{det}(\vect{\Sigma}_{\theta}^i)\big| + \frac{1}{2}(\vect{x}^i(t_{k+h}) - \vect{\mu}_{\theta}^i)^\text{T}\vect{(\Sigma}_{\theta}^i)^{-1}(\vect{x}^i(t_{k+h}) - \vect{\mu}_{\theta}^i)
\end{multline}
where the constant term is dropped since it does not affect the minimization. The training of the network is performed by minimizing the loss function $\mathcal{L}$ over the training set $\widetilde{\mathbf{D}}$.

\subsection{Order Parameters - Free Energy Landscapes}
\label{sec:order_parameters_energy}
Order parameters are coarse, collective variables that summarize the physics involved in the colloidal self-assembly process \cite{PhysRevB}. These quantities often encapsulate features of interest, for example the compactness, or the local or global degree of order of a particle assembly \cite{PhysRevB,tang2017construction}. Such variables can then be used to formulate (ideally analytical, but practically, here, data-driven) models to study the collective dynamics of complex systems. Domain scientists often have prior knowledge of good candidate order parameters based on experience, intuition, or mathematical derivations, and validate a good variable choice among such candidates \cite{tang2017construction}. 
Order parameters that are typically used to study colloidal self-assembly are $Rg$, $\psi_6$ and $C_6$ \cite{edwards2014colloidal}. $Rg$ is the radius of gyration, which quantifies whether the particle ensemble is expanded in a fluid state or condensed in a crystalline state; $\psi_6$ is the degree of global six-fold bond orientation order (near 0 for ideal gas, 1 for perfect single domain crystal); $C_6$ is the ensemble average of a local order parameter based on the number of neighbors each particle has with 6-fold order (0 for no neighbors with 6 fold-order, 1 for 6 neighbors with 6-fold order). Expressions for each of these are included in our prior publications. \cite{edwards2014colloidal}.
In our case, we do not use these theoretical ``usual candidate" colloidal order parameters in our model construction. We allowed the data to determine which and how many collective variables are needed by using the Diffusion Maps scheme. We then attempt to establish explainability of our data-driven variables in terms of the theoretical ones $Rg,\psi_6, C_6$, see Section \ref{sec:latent_observables}.

The effective potential  $G(\vect{x})$ quantifies the free energy landscape. It is obtained from the equilibrium probability distribution, the steady-state solution of the Fokker Planck equation \cite{risken1996fokker}. The integral equation for this effective potential (alternatively, potential of mean force or effective free energy), is given (up to a constant) by the equation \cite{risken1996fokker}
\begin{equation}
\label{eq:eff_potential}
    G(\vect{x}^i)  = -kT \int_{0}^{\vect{x}^i} \bigg[2(\mathbf{\sigma}^{2})^{-1}({\vect{\nu}} - \nabla \cdot \frac{\mathbf{\sigma}^2}{2})  \bigg] \cdot d\vect{r}.
\end{equation}
where $\vect{\nu}$ is drift and $\mathbf{\sigma}$ is the diffusivity. For our computations we chose the origin as the reference state.

\section{Results}

\subsection{Latent Observables}
\label{sec:latent_observables}
We start by coarse-graining Brownian Dynamics simulation (details about the simulations and the sampling are in the Appendix). 
Diffusion Maps discovers two latent non-harmonic coordinates denoted as $\phi_1$,$\phi_2$ \cite{dsilva2015parsimonious}. This suggests that two Diffusion Maps coordinates are enough to provide a more parsimonious representation of the original data set. The selection of those two Diffusion Maps coordinates was made by applying the local linear regression algorithm suggested by the authors in \cite{dsilva2015parsimonious}. We first check the interpretabilty of these data-driven observables by coloring the Diffusion Maps coordinates as functions of the three order parameters $Rg$, $\psi_6$ and, $C_6$. Those order parameters are physically meaningful coarse variables that measure the degree of condensation of the material (Section \ref{sec:order_parameters_energy}). 
It is worth highlighting that no pair of the three order parameters is exactly one-to-one with the Diffusion Maps coordinates;  however a clear trend appears: condensed configurations arise at the center of our manifold embedding, while further out from the center  disordered/fluid like structures are observed. This implies that our latent coordinates encode the physics of the Brownian Dynamic Simulations.
\begin{figure}[ht]
\begin{center}
\includegraphics[scale=0.35]{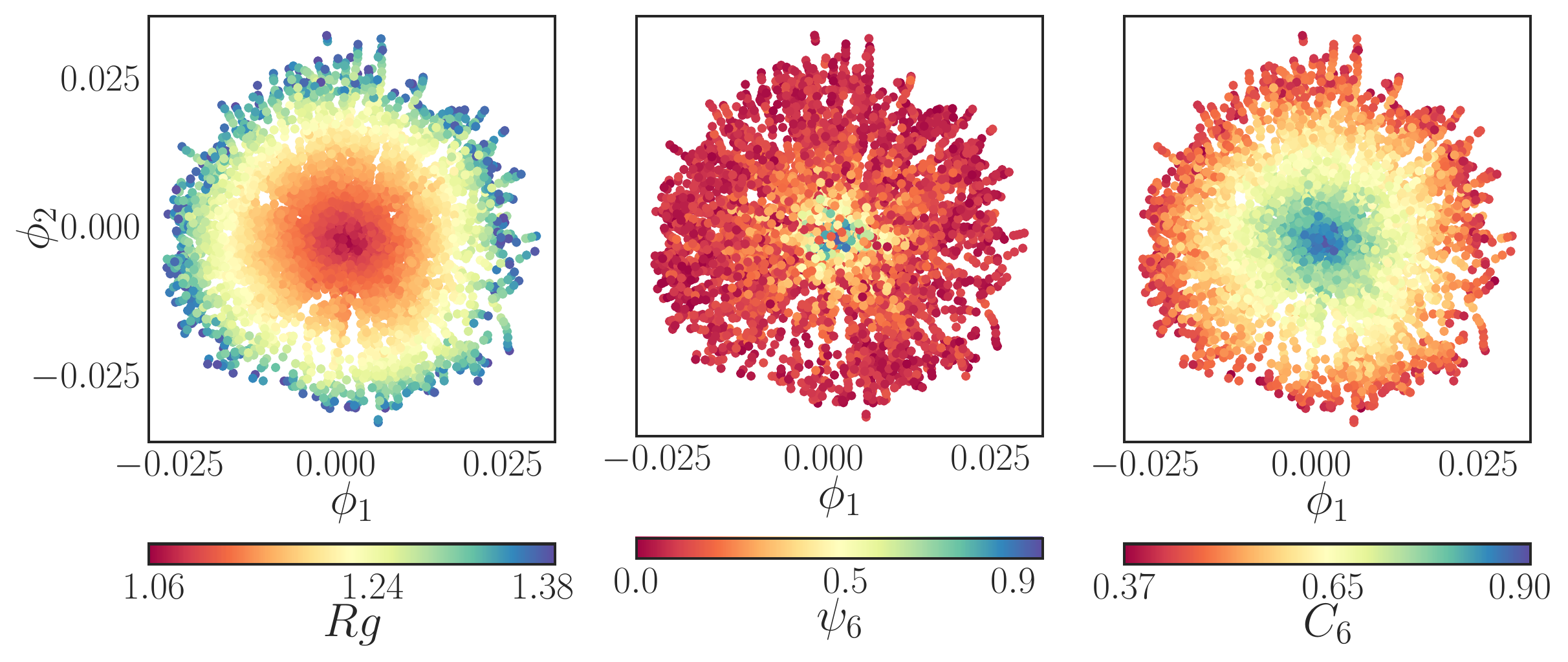}
\caption{The two leading Diffusion Maps coordinates ($\phi_1,\phi_2$)} colored by the physically meaningful order parameters $Rg$, $\psi_6$ and, $C_6$ respectively ($Rg$ clearly correlates with $C_6$). 
\label{fig:Diffusion_Maps_Order_Parameters}
\end{center}
\end{figure}

\subsection{Learning Effective SDEs}
\label{sec:learning_sdes}
\begin{figure}[ht]
\begin{center}
\centering
\includegraphics[scale=0.38]{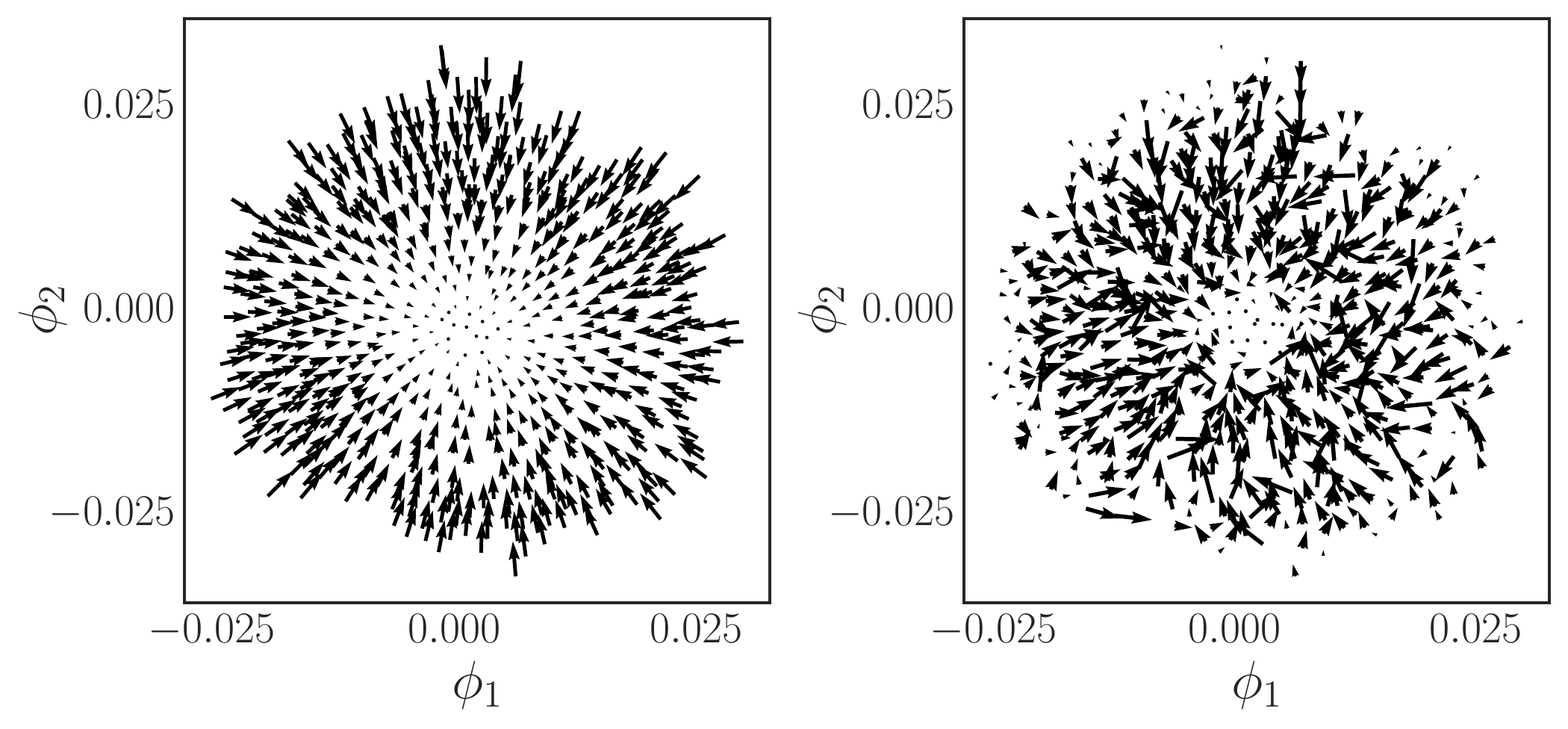}
\caption{The estimated drift from the neural network and the Kramers-Moyal respectively is plotted as the on average vector field in the two Diffusion Maps coordinates ($\phi_1,\phi_2$). The vector field is plotted in subsampled data sets to improve visualization.}
\label{fig:Vector_Field_Comparison}
\end{center}
\end{figure}

For Brownian Dynamics simulations at fixed normalized voltage \hbox{$V^* = \frac{V}{V_{xtal}} = \frac{1.51 \text{V}}{1.89\text{V}} = 0.8$} (see Section A2 in the Appendix), we estimate the drift and diffusivity in the Diffusion Maps coordinates with (a) a neural network architecture; and with (b) the Kramers-Moyal expansion. 
The drift estimated by the two approaches is plotted as a vector field on the two Diffusion Maps coordinates, Figure \ref{fig:Vector_Field_Comparison}. 
The drift component gives us an  estimate of what the trajectories will locally tend to do on average.
As can be seen from Figure \ref{fig:Vector_Field_Comparison} (on average) the trajectories will evolve towards the center of the manifold, and therefore towards more condensed structures as expected from the detailed Brownian Dynamics simulations.

\begin{figure}[ht]
\begin{center}
\centering
\includegraphics[scale=0.45]{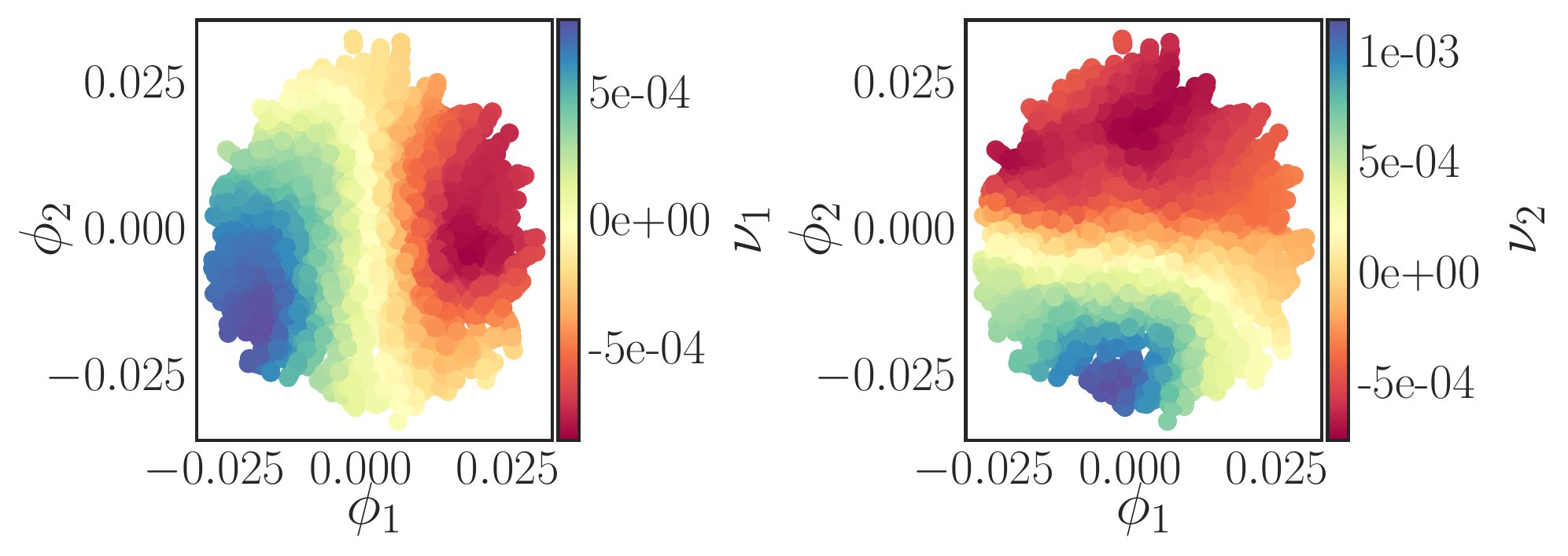}
\includegraphics[scale=0.45]{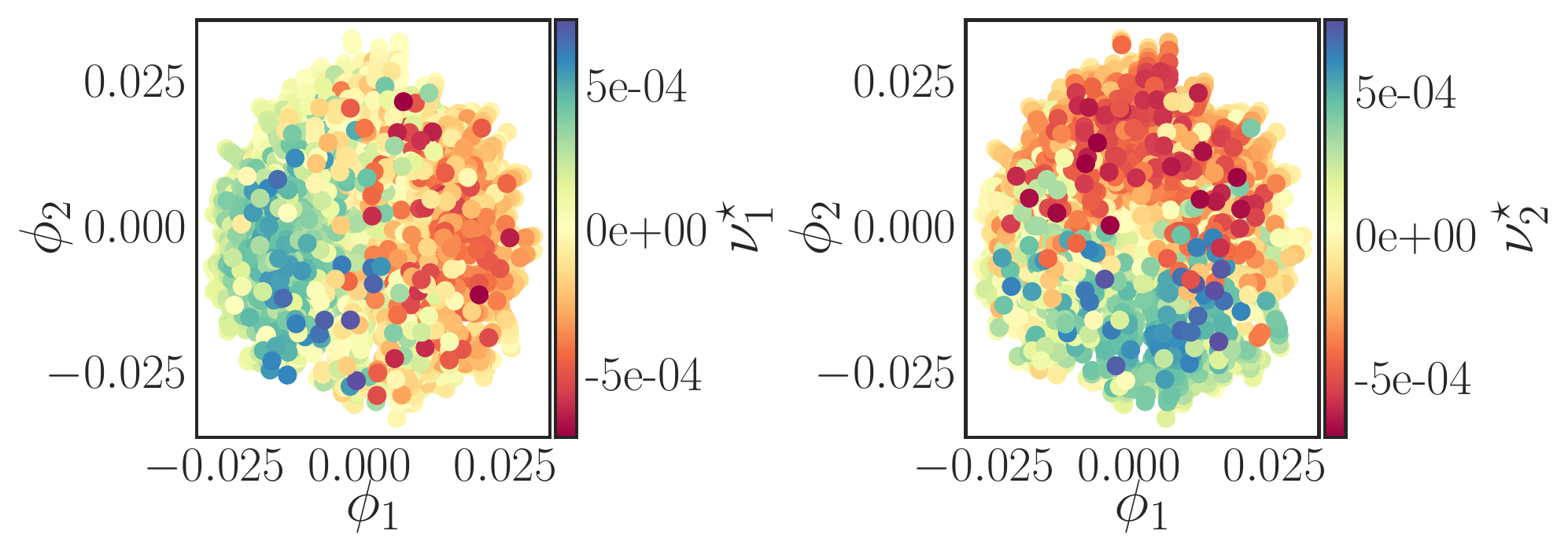}
\caption{The estimated drift, $\nu_1,\nu_2$, from the neural network (first row) and the Kramers-Moyal, $\nu_1^{\star},\nu_2^{\star}$, (second row) is plotted as a function of the Diffusion Maps coordinates ($\phi_1,\phi_2$)}.
\label{fig:Drift_Comparison}
\end{center}
\end{figure}

\begin{figure}[ht]
\begin{center}
\centering
\includegraphics[scale=0.35]{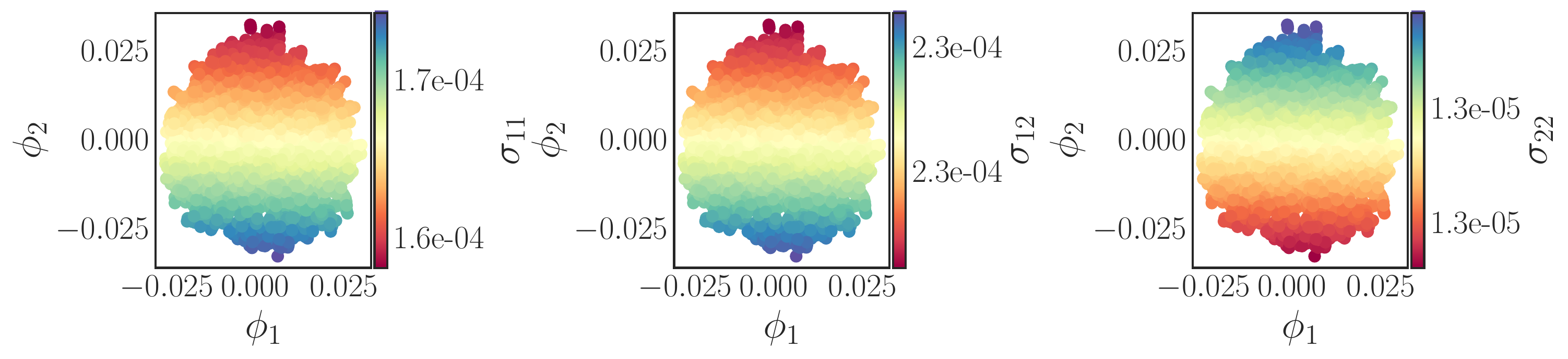}
\includegraphics[scale=0.35]{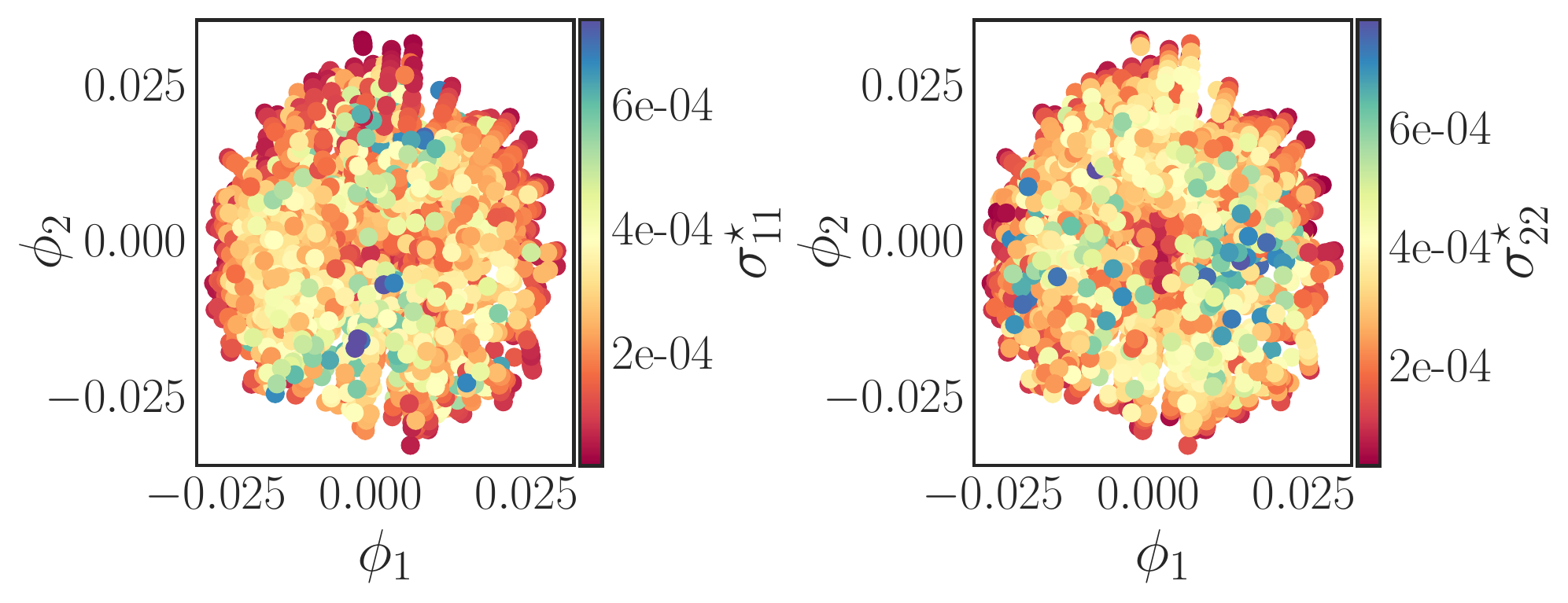}
\caption{The estimated diffusivity from the neural network, $\sigma_{11},\sigma_{22},\sigma_{12}=\sigma_{21}$ (first row) and the one from Kramers-Moyal, $\sigma^{\star}_{11},\sigma^{\star}_{22}$ (second row) is plotted as a function of the Diffusion Maps coordinates ($\phi_1,\phi_2$). } 
\label{fig:Diffusivity_Comparison}
\end{center}
\end{figure}

The estimated vector field from the neural network appears \textit{smoother} compared to the one obtained with the \hbox{Kramers-Moyal}. This could be partially attributed to the fact that the neural network during training for the drift learns simultaneously from many points per iteration through the loss function. On the other hand, Kramers-Moyal uses bursts around each individual data point separately, without information about the nearby points. 
Figure \ref{fig:Drift_Comparison} offers another comparison between the estimated drift from the neural network and the Kramers-Moyal Expansion. The estimated drifts of the two methods are comparable, with the drift estimated from the neural network often slightly larger in magnitude. 
The comparison for the estimated diffusivity with the two approaches leads to similar conclusions, Figure \ref{fig:Diffusivity_Comparison}. The diffusivity estimated from the neural network appears smoother compared to the one estimated from the Kramers-Moyal. 
Note the neural network approach estimated the diffusivity matrix without assuming it to be diagonal (as opposed to its Kramers-Moyal estimation). Even though a trend appears in the diffusivity computed through the network the computed diffusivity is practically constant along the data and the trend is just an artifact of the fitted diffusivity through the network.

Given the estimated drift and diffusivity we wish to generate trajectories for the reduced eSDEs in Diffusion Maps coordinates. Evaluating the drift and diffusivity along the integration is trivial for the trained neural network. For the Kramers-Moyal expansion, interpolating from the computed values becomes necessary. Since the functions of the estimated drift and diffusivity are not smooth enough for a global interpolation scheme, a local nearest neighbor interpolation was used during the integration. The numerical integrator used for both cases was the Euler-Maryama scheme.

From the estimated coefficients (drift and diffusivity), and more precisely from the average vector field in Figure \ref{fig:Vector_Field_Comparison}, it is expected that the trajectories will evolve toward the center of the embedding for both estimated eSDEs.

\begin{figure}[ht]
\centering
\begin{center}
\includegraphics[scale=0.35]{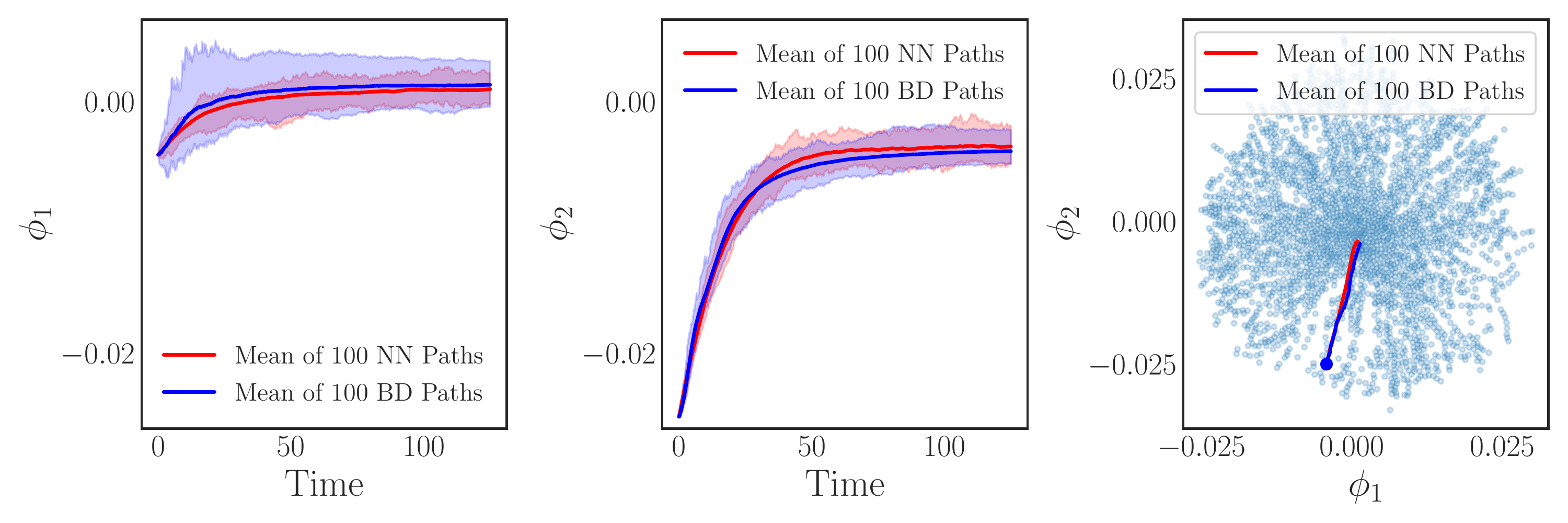}
\includegraphics[scale=0.35 ]{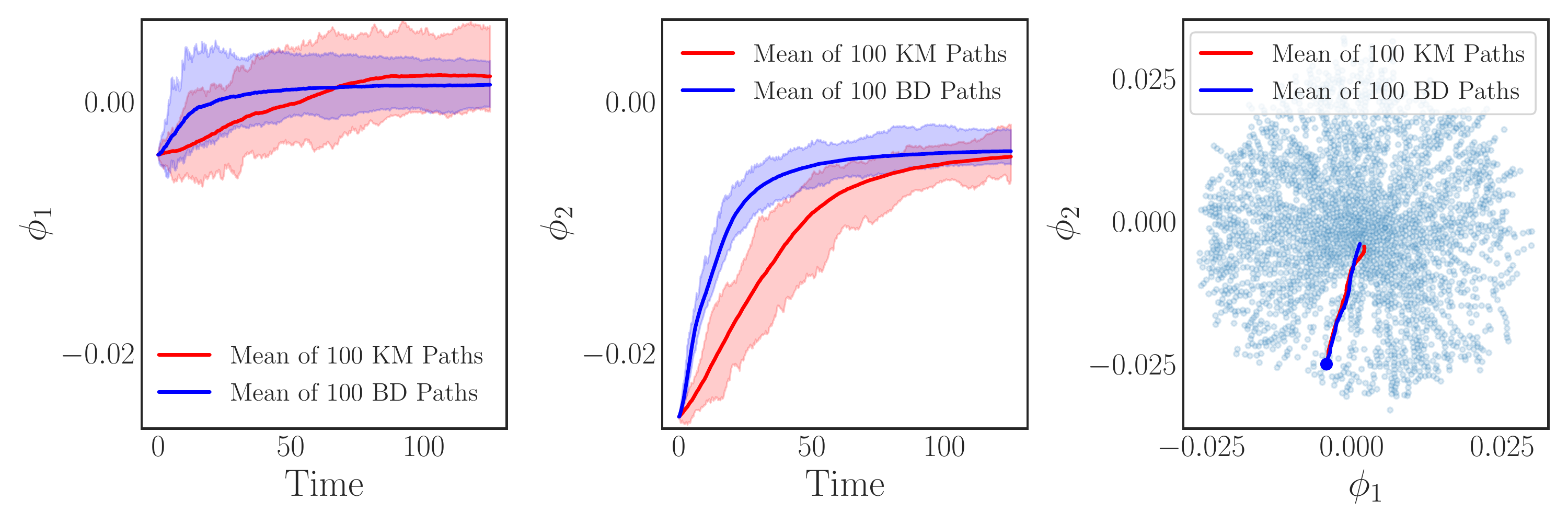}
\caption{The estimated dynamics from the neural network (first row) and the Kramers-Moyal (second row) is shown compared to restricted (with Nystr\"om) trajectories of the Brownian Dynamics in the two Diffusion Maps coordinates ($\phi_1,\phi_2$). The mean of 100 trajectories (starting from the same initial condition) is used for all cases. 
To get a visual inspection of the variance in the estimated 100 trajectories, the area between the maximum and minimum values for those trajectories is being ``filled'' with solid color.
The red paths correspond to the data-driven eSDEs (neural network or Kramers-Moyal) and the blue paths to the restricted Brownian Dynamics Paths.}
\label{fig:KM_Network_Comparison}
\end{center}
\end{figure}

To evaluate our models' performance against ground-truth data, we sampled Brownian Dynamics trajectories and \textit{restricted} those trajectories with Nystr\"om Extension in the reduced Diffusion Maps coordinates ($\phi_1,\phi_2$). A comparison between the mean of 100 trajectories obtained from the two eSDEs is contrasted to the mean of 100 restricted trajectories computed with Brownian Dynamics simulations in \hbox{Figure \ref{fig:KM_Network_Comparison}}. The dynamics from the network on average provide more accurate results compared to the ones obtained from the Kramers-Moyal.
To successfully estimate both the drift and the diffusivity for the neural network, training was performed in two stages. First, we chose a time step $h$ that gave a reasonable estimation of the drift;  we then \textit{fixed} the part of the network that estimates the drift, and used  snapshots at smaller time steps $h'$ to estimate the diffusivity (see the discussion in \cite{dietrich2021learningsdes}).

We provide an uncertainty quantification (error analysis) comparison of the neural network model  in Section A6 of the SI. This analysis provides some more quantitative measurements on of how certain the reported predictions are. The results suggest the the robustness of the neural-network model. In addition, in Section A7 of the SI we discuss a more quantitative comparison between the two surrogate identified eSDEs (with Kramers-Moyal and the neural-network). The results in this case suggest that the discrepancy between the two surrogate models is in the range of the expected error estimations of the neural-network model.

\subsection{Learning a Parameter-Dependent eSDE}
\begin{figure}[ht]
\centering
\includegraphics[scale=0.35]{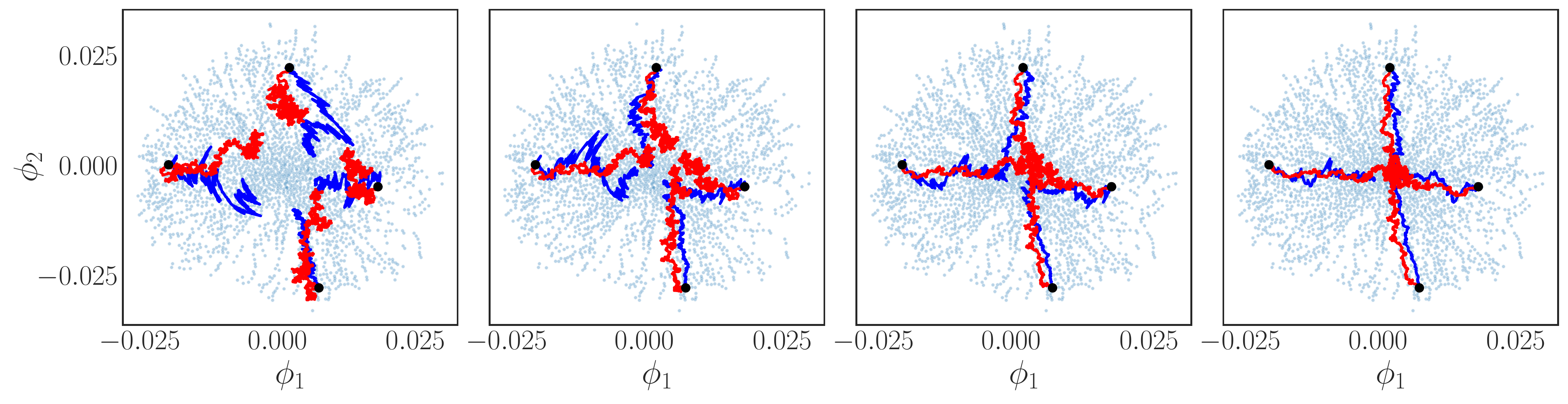}
\caption{In the Diffusion Maps coordinates ($\phi_1,\phi_2$) we illustrate trajectories computed through the neural network trained for different values of the voltage, $V^{\star} = \{0.5,0.6,0.7,0.8 \}$ plotted from left to right. Trajectories of the estimated eSDE from the neural network (red paths) are contrasted to restricted (with Nystr\"om Extension) trajectories computed with Brownian Dynamics (blue paths) for different values of the voltage.} 
\label{fig:Phase_Portait_Four_Voltages}
\end{figure}

In this section we illustrate the ability to learn a parameter dependent eSDE. For this case only the neural network was used. We sampled data (snapshots) for four different voltages, \hbox{$\mathbf{V}^{\star} = \{0.5,0.6,0.7,0.8\}$}.
The larger the voltage becomes, the larger the force that is acting on the particles, and thus the faster they condense. On the contrary, as the voltage becomes lower, the particles can move more freely and they condense slower. Those physical features are expected to be captured in terms of the drift and diffusivity of our eSDE . As the voltage increases the drift (force) is expected to increase and the diffusivity to decrease. In Figure \ref{fig:Phase_Portait_Four_Voltages} the obtained results from the neural network appear to conform to those features of the simulations. 

For four different initial conditions, and for the same time length, trajectories were integrated with Euler-Maryama; the same integration step was used for all parameter values. As the parameter value increases, from left to right in Figure \ref{fig:Phase_Portait_Four_Voltages}, the trajectories appear to travel faster towards the center of the embedding (towards more condense configurations). This can be attributed to fact that the drift increases in magnitude. In addition, as the voltage decreases, the trajectories appear more noisy,  since the diffusivity increases. 

\begin{figure}[ht]
\vskip 0.2in
\begin{center}
\centering
\includegraphics[scale=0.35]{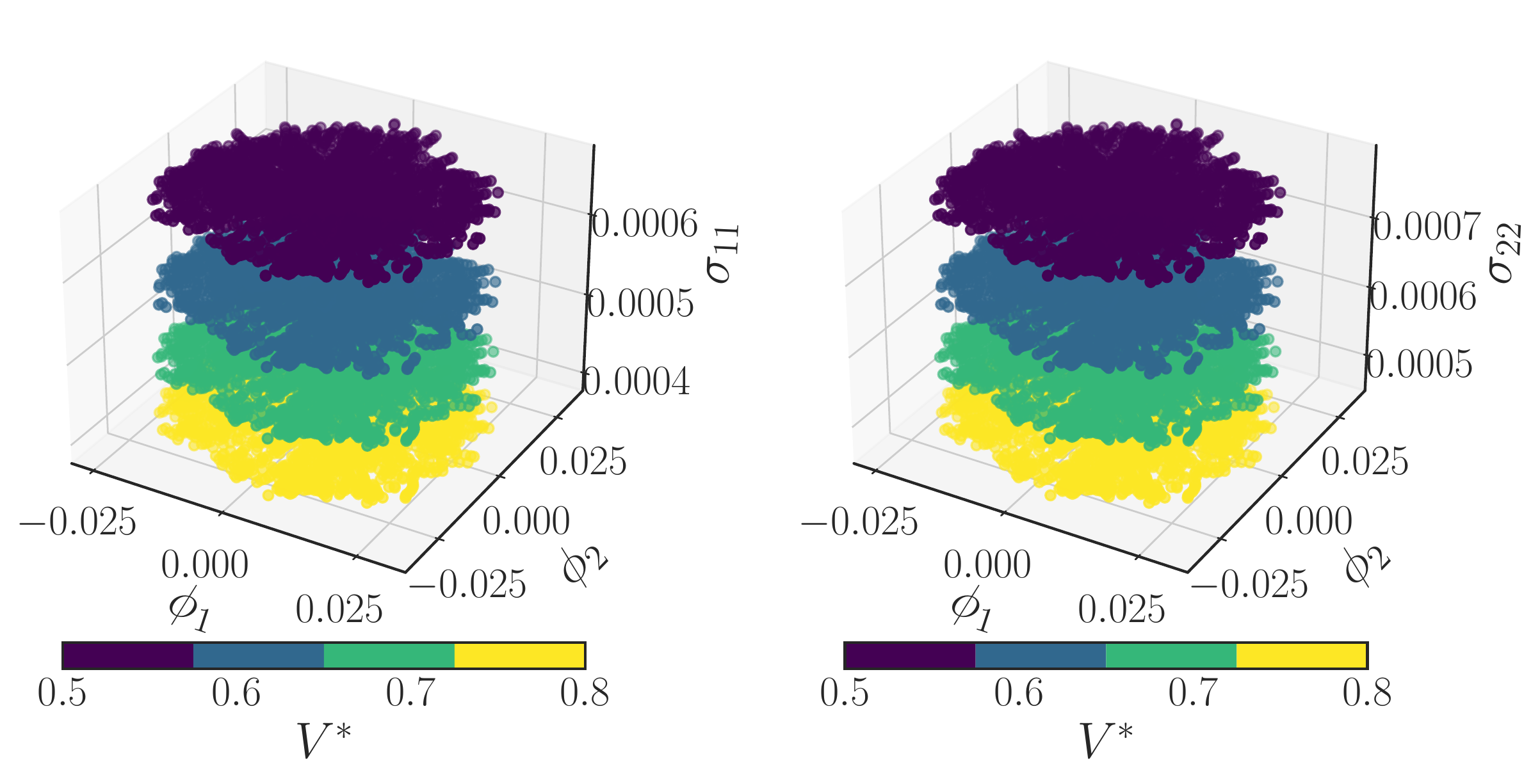}
\caption{The estimated diagonal diffusivity is plotted against $\phi_1$, $\phi_2$. Different colors correspond to different values of the parameter  ($V^*$). The trend in the estimated diffusivity is in agreement with the physics of the problem. The higher voltage forces the system to condense faster. On the contrary, the smaller the voltage becomes, the easier it is for the particles to move freely, and thus the effective diffusivity increases.} 
\label{fig:Diffusivity_Four_Voltages}
\end{center}

\vskip -0.2in

\end{figure}

In Figure \ref{fig:Diffusivity_Four_Voltages} the estimated diffusivity is plotted against the Diffusion Maps coordinates and is colored with the voltage value. Figure \ref{fig:Diffusivity_Four_Voltages} supports the observation that as the voltage increases the diffusivity decreases.

For the estimation of the parameter dependent eSDE, the flexibility of having different step sizes $h^i$ proved quite useful. For smaller values of the voltage $V^{\star}$, for which the drift is also smaller, larger time steps could be accurately employed.

\subsection{Free Energy Landscapes}
\begin{figure}[h]
\centering
\begin{center}
\includegraphics[scale=0.35]{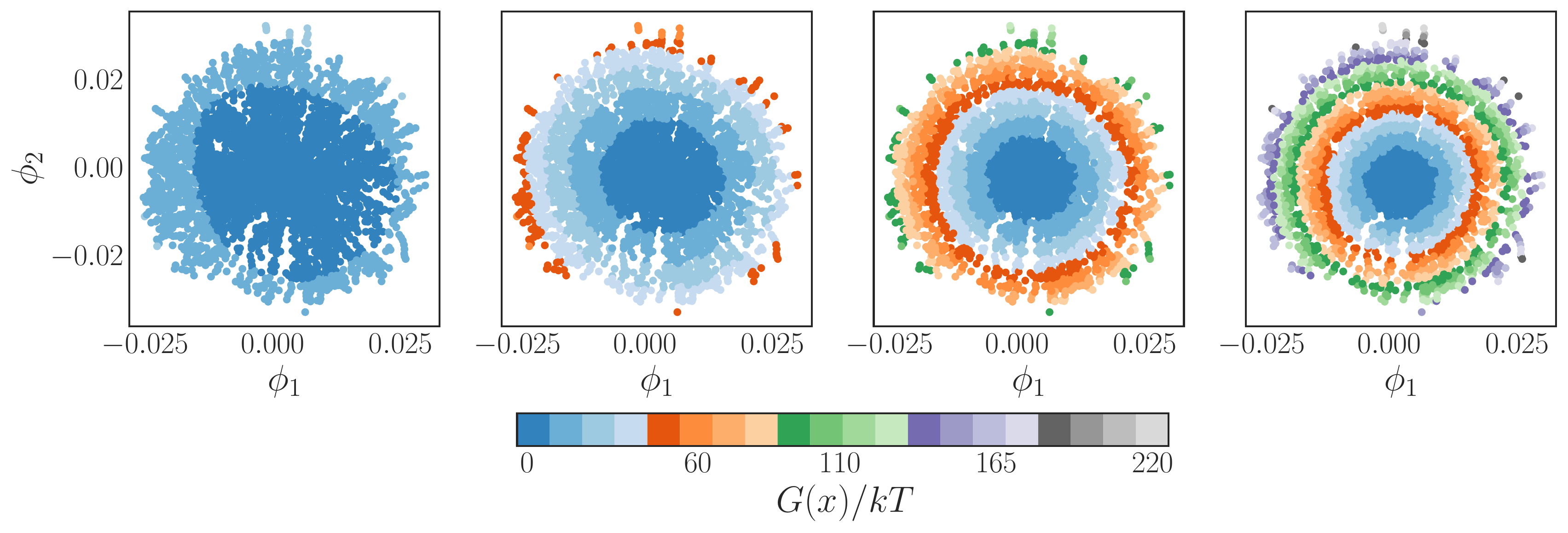}
\caption{The Free Energy Landscape, $G(x)/kT$, estimated by Equation \eqref{eq:eff_potential} for different voltages $\vect{V}^{\star} = \{0.5, 0.6, 0.7, 0.8 \}$ (from left to right) are plotted as functions in the Diffusion Maps coordinates ($\phi_1,\phi_2$) respectively.} 
\label{fig:Effective}
\end{center}
\end{figure}

We illustrate the ability to estimate free energy landscapes (potential functions) from the coefficients of the reduced eSDE. 
In Figure \ref{fig:Effective} the Free Energy in $kT$ units is plotted as a function of the Diffusion Maps coordinates $\phi_1$, $\phi_2$ for the four different voltages in an increasing order. From Figure \ref{fig:Effective} the larger the voltage, the larger the range of effective potential values becomes. From our computations it appears that the term $\nabla \cdot \mathbf{\sigma}^2$ is negligible compared to the other terms, and that the state dependence of the diffusivity can be practically ignored.

\subsection{Experimental Data}
In this section we provide a qualitative comparison between our reduced model and experimental dynamic data. The experimental set up from which the data were collected is described in \cite{edwards2013size}. Note that each configuration used for the experimental data has 204 particles and not 210 as in our simulations. The radii of the particles is the same as the one used for the simulation and the voltage ($V^{\star} = 0.74)$ for those experiments is in the range  of the voltages used to train the parameter dependent eSDE. Given the experimental trajectories, we used the same preprocessing as for the computational data, and then Nystr\"om Extension was used to restrict the experimental  configurations in the Diffusion Maps coordinates. The experimental data were rescaled in the same range as the simulations based on the ratio of the radii of the two reference configurations used for the Kabsch algorithm. Please note that to restrict the configurations with 204 particles in the Diffusion Maps coordinates obtained from the simulations, a different reference configuration was used for the Kabsch algorithm. The reference configuration was selected also here as the configuration with the smallest value of $Rg$ from the experimental trajectories. Then the same density estimation described in Section \ref{sec:Diffusion_Maps} is applied. These steps allow us to project the experimental particles to the Diffusion Maps coordinates despite their different number of particles. We then use our trained neural network to generate trajectories given the estimated initial conditions in the Diffusion Maps coordinates. The integration of the eSDE was performed for 125 seconds with time step $h =0.125$. The time step used to integrate the eSDE corresponds to the same frame rate that the experimental measurements were sampled at (8 frames per second \cite{edwards2013size}). The behavior of the restricted experimental trajectories has the same qualitative behavior as the reduced model, and as the restricted trajectories of the Brownian Dynamics.
\begin{figure}[ht]
\begin{center}
\includegraphics[scale=0.35]{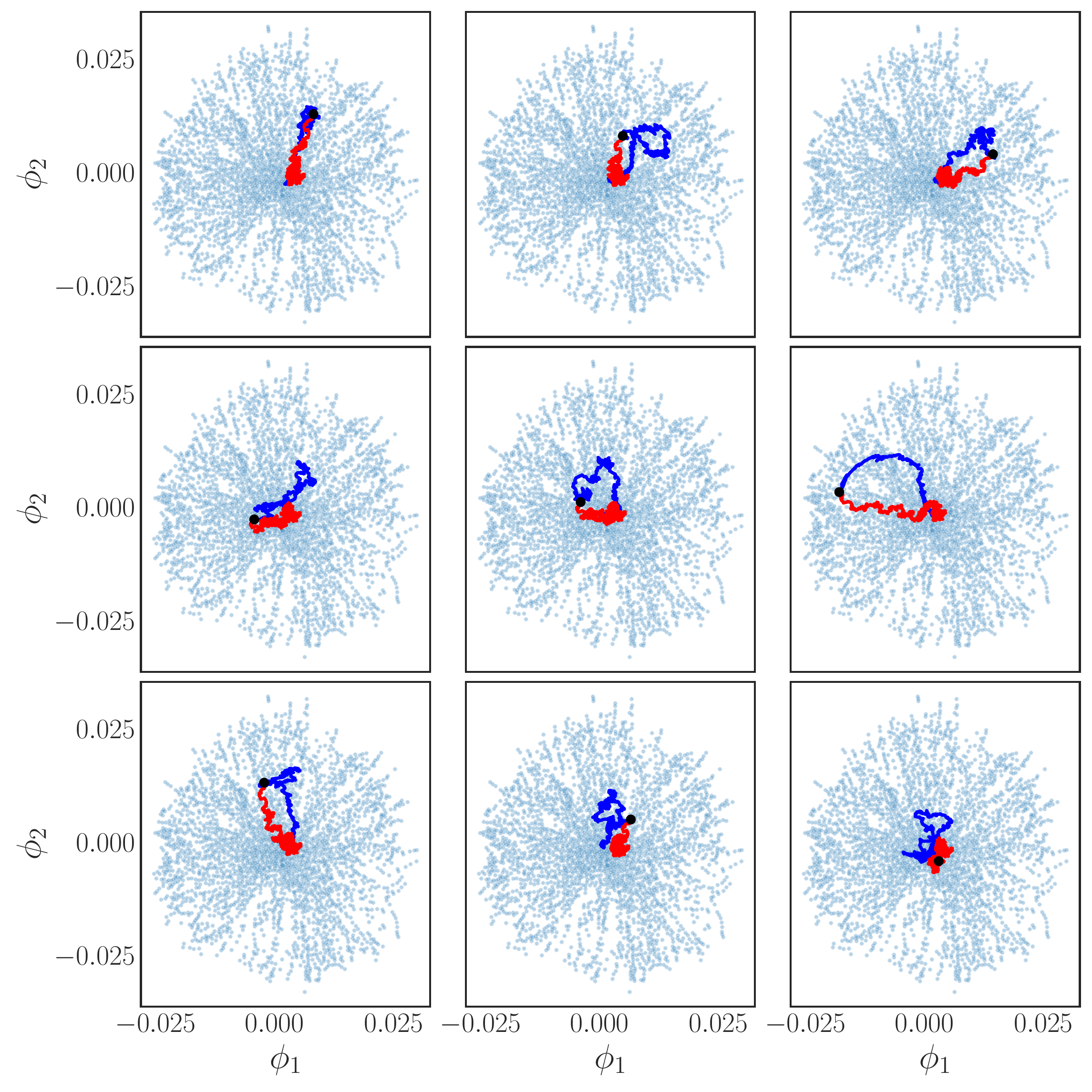}
\caption{Restricted with Nystr\"om Extension experimental trajectories compared with paths generated from the neural network eSDE trained on the computational data. The red trajectories correspond to the paths generated by the neural network and the blue trajectories to the trajectories of the experiments, restricted to the latent space using Nystr\"om Extension.}
 
\label{fig:Experimental_Comparison}
\end{center}
\end{figure}

\section{Discussion}
We demonstrated that the Diffusion Maps algorithm can discover a set of latent observables given a data set of sampled dynamic configurations of crystallizing colloidal particles. We explored the correspondence between our obtained latent observables and established theoretical order parameters ($Rg$, $\psi_6$, $C_6$). 
We learned an eSDE by using the traditional Kramers-Moyal expansion and compared it with a modern deep learning architecture based on stochastic numerical integrators \cite{dietrich2021learningsdes}. Both estimated reduced eSDEs qualitatively reproduce the dynamics of the full simulations. We showed that the neural network's dynamics on average appear more accurate, by comparing the data-driven eSDEs with restricted trajectories of the full Brownian Dynamics. It's worth mentioning that the computation cost and the number of data points needed to learn an effective eSDE with the neural network is much smaller than the corresponding Kramers-Moyal effort, we provide a more detailed comparison in the Appendix. We illustrated the ability to learn a parameter dependent eSDE through our neural network architecture. The coefficients (drift and diffusivity) of the parameter dependent eSDE again seem to capture the dynamics of the fine scale simulations. Lastly, we showed that our reduced models qualitatively agree with dynamics of restricted experimental data.

\section{Conclusions}
The developed eSDEs provides a compressed data-driven model that we believe can help the study of self-assembly. Even though the application was focused on colloidal assembly, this framework can be applied to a range of different applications, from coarse-graining epidemiological models to models of cell motility. Such data-driven models could be useful tools for performing scientific computations (e.g. estimation of mean escape times, construction of bifurcation diagrams) even when analytical expressions are not available.

\begin{figure}[h]
\begin{center}
\centering
\includegraphics[scale=0.35]{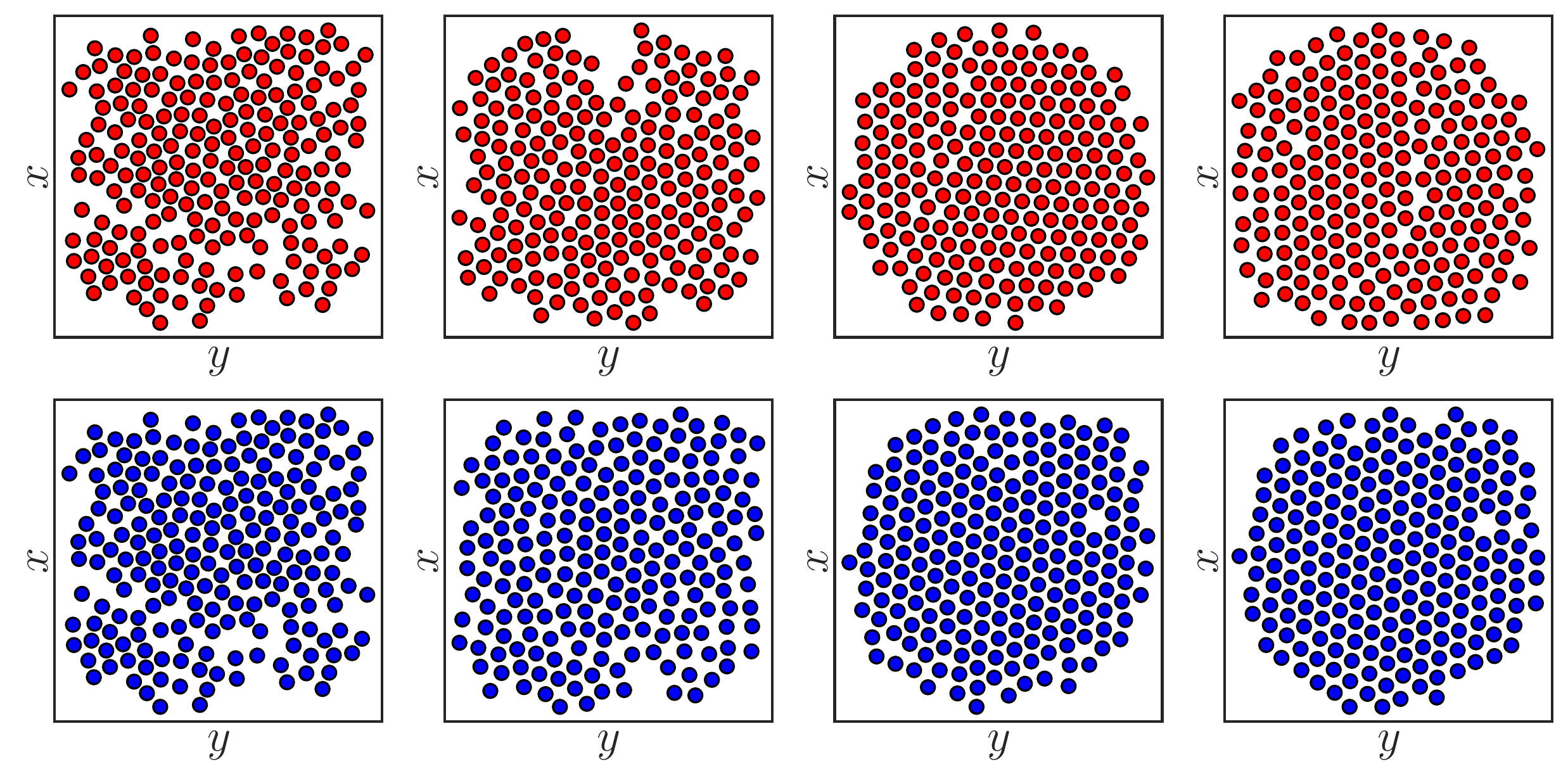}
\caption{The first row illustrates snapshots of the colloidal particles at different time instances ($t=0 \text{s}, 12.5 \text{s}, 62.5 \text{s}, 125 \text{s}$). Those snapshots generated by mapping a trajectory integrated by the eSDE to the original physical coordinates with the nearest neighbor algorithm (lifting).} The second row illustrates snapshots of colloidal particles for the same timestamps computed with Brownian Dynamics. 
\label{fig:Lifting}
\end{center}

\end{figure}

Our reduced models, while capable of describing the coarse-grained, collective dynamics, do not provide information about the fine-scale conformations themselves. 
Our assumption that differences in density profiles suffice to determine a similarity measure in configuration space leads to configurations with the same density field being mapped to a single point in our coarse latent space. Therefore, mapping back to the ambient space, i.e. \textit{lifting}, is a nontrivial task since there is a \textit{family} of configurations for each Diffusion Maps point. To support this argument we show a comparison between a \textit{naive} mapping of a generated trajectory from the Diffusion Maps coordinates to the configurations with nearest neighbors (what we call \textit{lifting}, from coarse to fine) and a trajectory generated by the Brownian Dynamic simulation
in Figure \ref{fig:Lifting} and in the accompanying video (provided in the SI). Both trajectories start from the same initial condition. In Figure \ref{fig:Lifting} for a reduced trajectory of the eSDE estimated by the neural network we find the nearest point in the Diffusion Maps coordinate that belongs to our data set and \textit{lift} based on that configuration. The \textit{lifted} trajectory exhibits large abrupt changes and appears thus unrealistic. We believe that utilizing conditional Generative Adversarial Networks (cGANs) constitutes a promising direction for reconstructing realistic fine scale configurations 
conditioned on coarse-grained features.

 Learning eSDEs directly from experimental data is a possible extension of our work. The main limitation of learning an eSDE directly is that usually we do not have sufficient experimental data; this is why BD models are matched (as well as possible) to experiments, and then we analyze their simulations\cite{edwards2013size,juarez2009interactions,juarez2011kt}. Perhaps a transfer learning approach, where the eSDE is initially trained in a large computational data set, and then refined/adapted to experimental data could be an interesting approach for the construction of data-driven models for studying self-assembly.

Another possible extension of our current work deals with using the identified eSDE for control problems. Merging the parameter dependent eSDE with feedback control policies could guide the evolution of configurations from polycrystalline states to target
    single-domain crystals \cite{zhang2020controlling,tang2016optimal}.
 
\section*{Author Contributions}
    N.E. performed the data curation, the formal analysis of the data and conducted the experiments (simulations and learning the eSDEs) . F.D, provided the codes for the neural networks made the modifications for the architecture to include a parameter and provided guidance for the training and evaluation of the neural network models. J.M.B-R, provided theoretical and practical insight for the computations of the effective potentials.  A.Y, M.A.B provided the Brownian Dynamic Simulation codes and codes for the computation of the Order Parameters. R.S, M.A.B  analyzed and provided data for the experimental paths. I.G.K.  planned and coordinated the work and was the PI in funding the effort.  N.E, F.D, and I.G.K wrote the manuscript. All authors contributed in manuscript's review and editing.

\section*{Conflicts of interest}
There are no conflicts to declare.

\section*{Acknowledgements}
This work was partially supported by the US AFOSR and the US Department of Energy. We acknowledge financial support by the National Science Foundation CBET-1928950.

\bibliographystyle{plain}
\bibliography{lit}

\newpage
\appendix
\onecolumn
\section{Appendix}
\subsection{Sampling}
\label{sec:Sampling}
The data set in which Diffusion Maps was computed was sampled as follows. Brownian Dynamic simulations were performed given a fixed voltage $V^{\star} = 0.8$ for a system of 210 particles. For $\sim$ 1500 \textit{random} initial configurations the system was integrated for 1000s. Particle configurations were stored every 1s. The order parameters $Rg$ and $\psi_6$ were computed based on these collected snapshots. The configurations with $Rg>1.39$ were discarded. Since these dilute states arise from the random initialization of each simulation, the remaining data set contains configurations that cover the fluid, polycrystalline and crystalline space. We further subsampled in $Rg,\psi_6$ space in order to get a more uniform data set. The subsampled data set in which Diffusion Maps was performed contains $\sim$11,000 configurations.

For the Kramers-Moyal expansion (discussed in Section 2.3.1) further subsampling of the data set used for Diffusion Maps was performed. A data set of $\sim 2800$ configurations (and subsequently Diffusion Maps coordinates) was used. For each point in this data set $300$ 
\textit{short} Brownian Dynamic simulations were performed.

For the deep learning approach (discussed in Section 2.3.2) the same data set as for the Kramers-Moyal expansion was used. In this case for each data point only a short trajectory was needed in order to compute the snapshots. The same data points were also used to sample snapshots for  different values of the parameter. Different time steps (h) were necessary in order to sample properly the snapshots for the different values of the parameters. The fact that the dimensionality of the manifold remained the same in this range of parameter values was checked before learning the parameter dependent eSDE.

The computational time needed to sample data points for the Kramers-Moyal compared to the one for the neural network is $\sim 38$ times and the number of data used to compute Kramers-Moyal was $\sim$ 100 times more than the one used for the neural network. There is no significant difference between the training time needed for the neural network and Kramers-Moyal estimation. 
\subsection{Parameters for simulations}
\label{sec:Parameters}
\begin{table}[th]
\caption{Simulation Parameters of Colloidal Particles in the presence of a Quadrupole Electrode.}
\label{sample-table}
\vskip 0.15in
\begin{center}
\begin{small}
\begin{tabular}{|p{7cm}|p{4cm}|}
\toprule
Variable & Value  \\
\midrule
Number of particles, $N$    & 210 \\
Particle Size, $\alpha$ (nm), & 1400 \\
Temperature, $T$ ($^\circ C$)  & 20   \\
Clausius-Mossotti factor for 1 MHz AC field, $f_{cm}$    & -0.4667\\
Debye length, $\kappa^{-1 }\text{nm}$     & 10 \\
Electrostatic potential prefactor, $B^{PP}$ (kT)      & 3216.5 \\
Lowest voltage to crystalize system, $V_{xtal}$ (V)      & 1.89  \\
Applied voltage, $V$(V)   & 0.95, 1.13, 1.32, 1.51 \\
Normalized voltage, $V^*$ & 0.5, 0.6, 0.7, 0.8\\
\bottomrule
\end{tabular}
\end{small}
\end{center}
\vskip -0.1in
\end{table}
\subsection{Video}
In the attached video in the upper left a trajectory that evolves by  integrating the eSDE is shown.
In the lower left, a trajectory integrated with the Brownian Dynamic Simulations and restricted with Nystr\"om Extension in the Diffusion Maps coordinates is shown. In the upper right the \textit{lifted} with nearest neighbors estimation configurations are shown. In the down left the configuration integrated with the Brownian Dynamic Simulations are shown.
\subsection{Codes}
All the codes that used to produce the Figures and generate the models are provided in the public repository  Gitlab repository \href{https://gitlab.com/nicolasevangelou/colloidal-assembly-learning-sdes.git}{here}. 
\subsection{Neural Network Models Architecture}
The architecture and the training protocols for our two neural-network models are reported in this Section. 
For the surrogate model with fixed voltage, $4$ layers with $25$ neurons were used. The activation function for the drift network was selected as \textit{ReLU} (activation function \textit{ELU} was also tested and gave similar results). The activation function for the diffusivity network was \textit{softplus}. The training of the model was made in two stages: first we use snapshots with $h=1.0$ and train for 200 epochs with batch size equal to 32. This first stage approximates reasonably well the drift, but gives a much larger diffusivity. To decrease the diffusivity estimate we fix the weights for the drift network to \textit{non-trainable}. We then use snapshots of $h=0.125$ and train for $1,000$ more epochs.

For the parameter dependent eSDE a network with $5$ layers and $26$ with \textit{ELU} was trained. This network was initially trained for $1,000$ epochs with batch size equal to 51. We then set the weights for the diffusivity part of the network to \textit{non-trainable} and train for the drift for $5,000$ epochs. We then set the weights for the drift part of the network to \textit{non-trainable} and train for $1,000$ epochs to improve the diffusivity estimations. For each voltage snapshots at different step size $h$ were used.

For both networks outliers were removed (e.g. with \textit{z-score}). This led to a better approximation of the diffusivity. To train each model, the data was split 90$|$10 (train/validation). As test set we used the individual trajectories reported in the main text.

\subsection{Uncertainty Quantification of our Neural-Network Model}
\label{sec:Uncertainty_Quantification_NN}
We performed the following computations to obtain an estimate of the accuracy and the sensitivity of the learned parameters of the neural network. We trained the neural network model 200 times, each time using different splits of our original data set to training/validation  sets. For each of those models in tensorflow we used \textit{validation\_split=0.1}. The data were shuffled before training. To ensure that each time the training/validation sets we get are unique, a different seed was used from the \textit{numpy} random generator for each of the 200 models. To alleviate the computational time needed to train these 200 networks we used \textit{job-arrays}. All the other hyperparameters involved in training the model were kept fixed during this procedure. To obtain a quantitative measurement of the sensitivity of our model in terms of the training set, we then evaluated the drift and diffusivity for the entire data set (training and validation) as well as on a grid of test points.

\begin{figure}
\begin{center}
\includegraphics[scale=0.4]{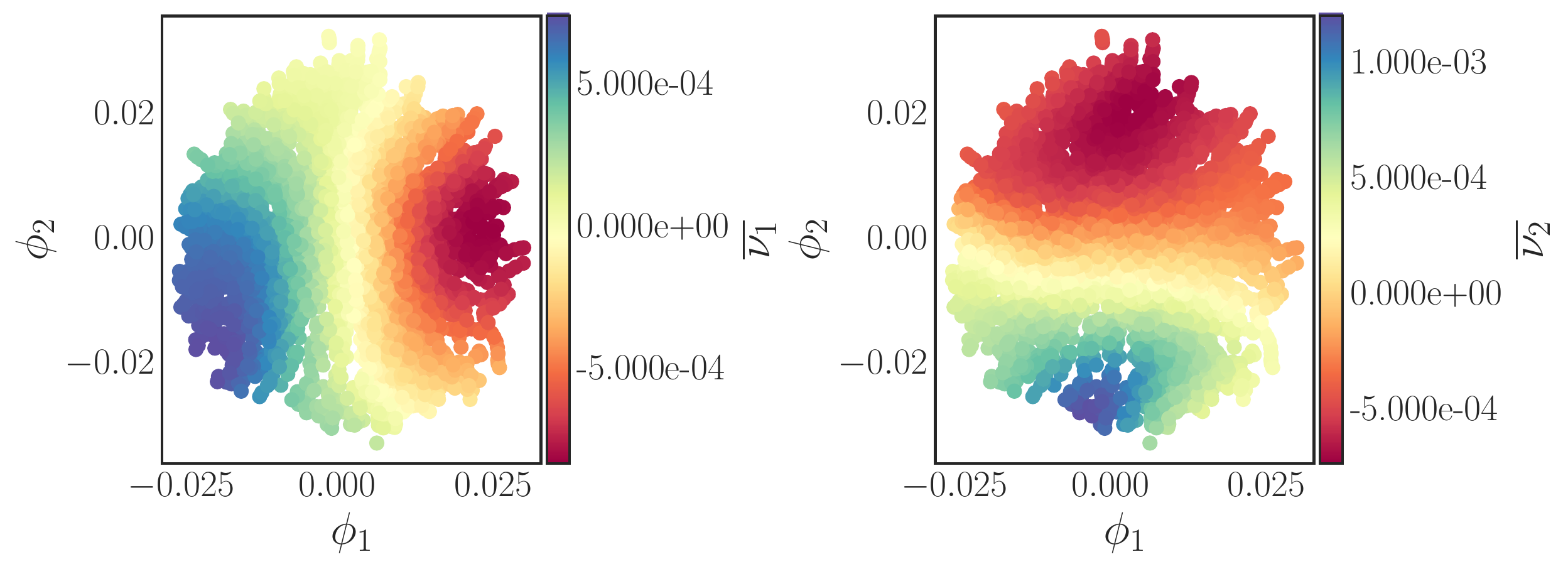}
\includegraphics[scale=0.4]{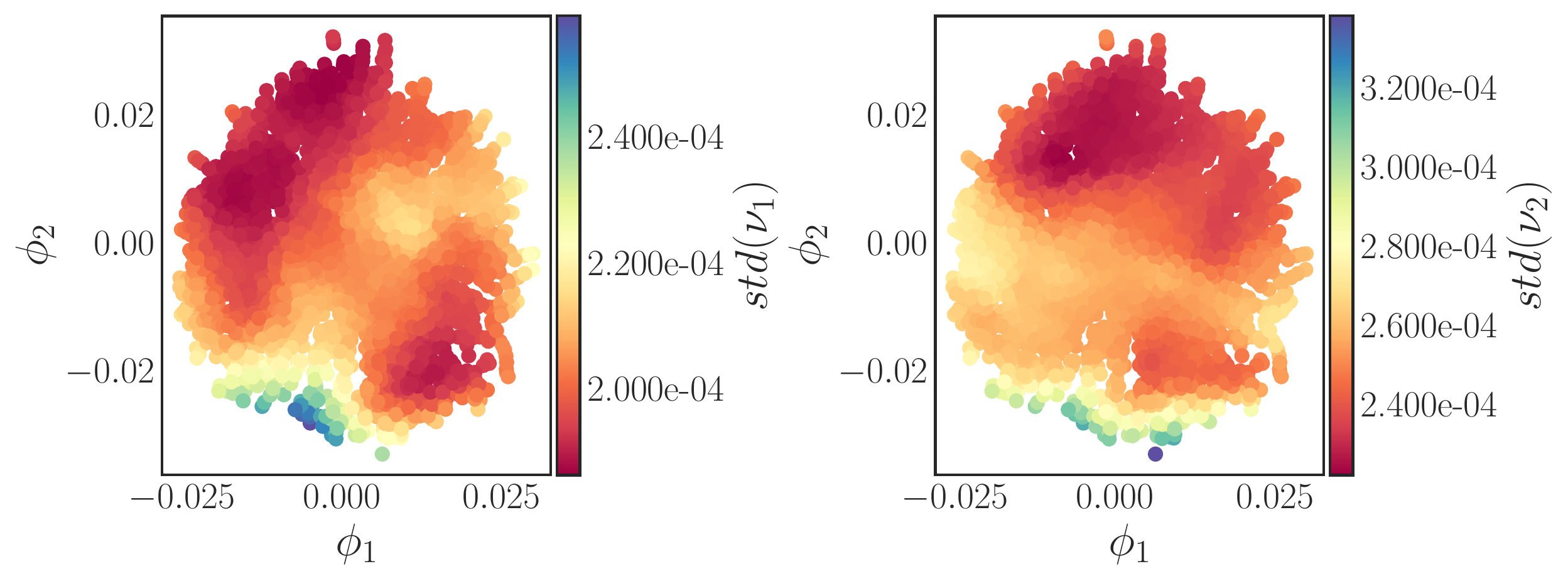}
\caption{[First row] The average values (over the 200 models) of the drifts $\overline{\nu_1}$,$\overline{\nu_2}$ are shown as functions in the Diffusion Maps coordinates. [Second Row] The standard deviation of the estimated drifts $std(\nu_1)$, $std(\nu_2)$ is shown as functions in the Diffusion Maps coordinates.}
\label{fig:drift_avg_dmaps}
\end{center}
\end{figure}

In Figure \ref{fig:drift_avg_dmaps} the first row illustrated the drift components (averaged over the 200 models) colored as a function of the diffusion maps coordinates. The second row depicts the pointwise standard deviation of the 200 models, colored as a function of the Diffusion Maps coordinates. It appears that the overall trend and the order of magnitude between the estimated average value of the drift here, and the one shown in Figure 4 in the main text are consistent. The estimated standard deviation appears almost everywhere smaller than the estimated mean drift.

\begin{figure}[ht]
\begin{center}
\includegraphics[scale=0.35]{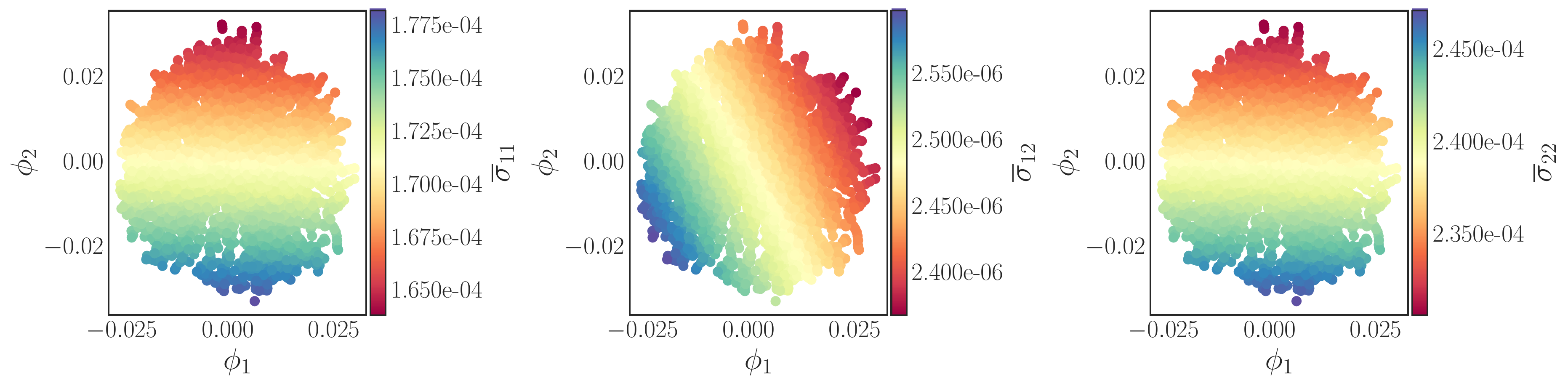}
\includegraphics[scale=0.35]{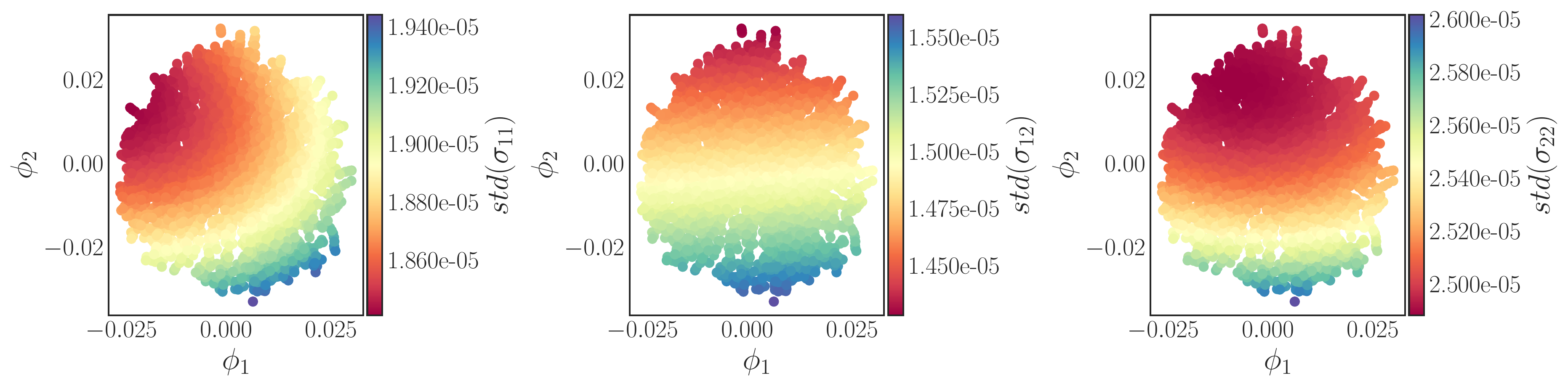}
\caption{[First row] The average values (over the 200 models) of the estimated diffusivities $\overline{\sigma_1}$,$\overline{\sigma_2}$ are shown as functions in the Diffusion Maps coordinates. [Second Row] The standard deviation of the estimated diffusivities $std(\sigma_1)$, $std(\sigma_2)$ is colored as a function in the Diffusion Maps coordinates.}
\label{fig:diffusivity_avg}
\end{center}
\end{figure}

The same computation was performed for the diffusivity and is shown in Figure \ref{fig:diffusivity_avg}. The magnitude of the diagonal terms is comparable to the model described in the main text (Figure 5). The average off-diagonal term seems to be an order of magnitude smaller than the one presented in the main paper. The average off-diagonal terms are two orders of magnitudes smaller than the average diagonal terms.
The standard deviations for the diffusivities for the diagonal elements are one order of magnitude smaller than the averages; yet for the off-diagonal elements, they are an order of magnitude larger than the average. This latter observation of the off-diagonal elements can be attributed to the fact that some models are having positive and some negative off-diagonal elements. Therefore their mean is closer to zero but their standard deviation is larger than zero.

\begin{figure}[ht]
\begin{center}
\includegraphics[scale=0.35]{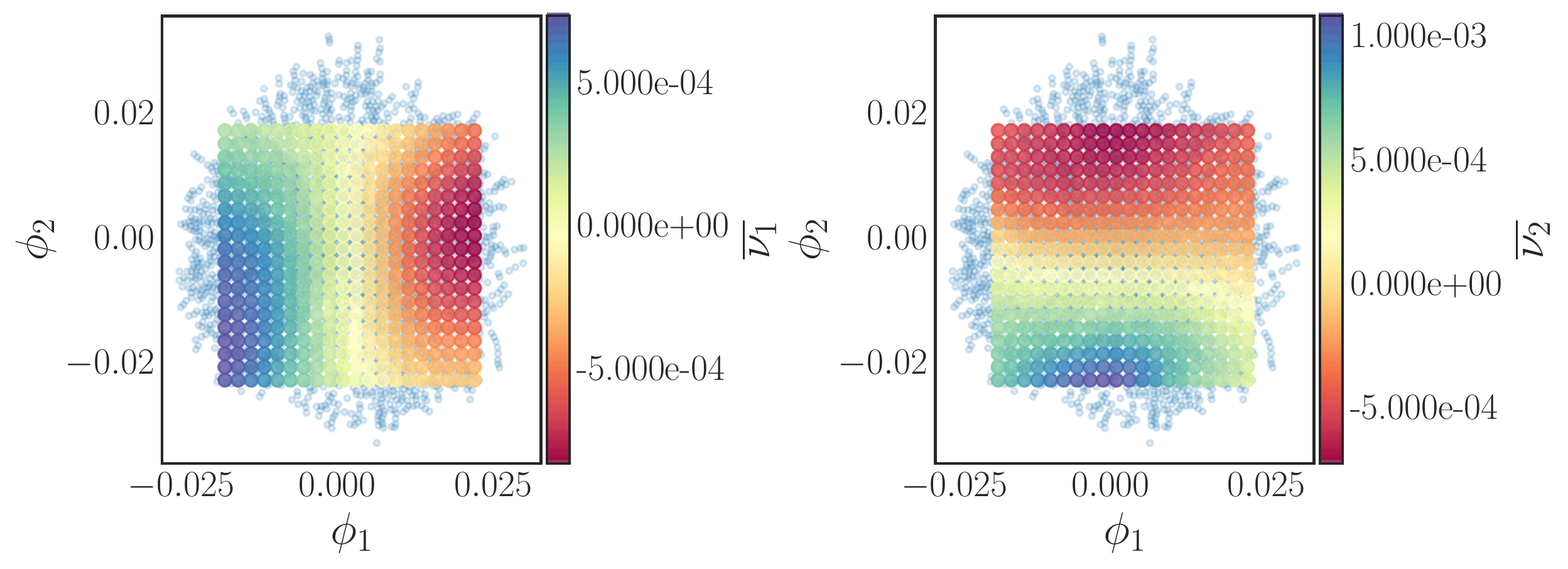}
\includegraphics[scale=0.35]{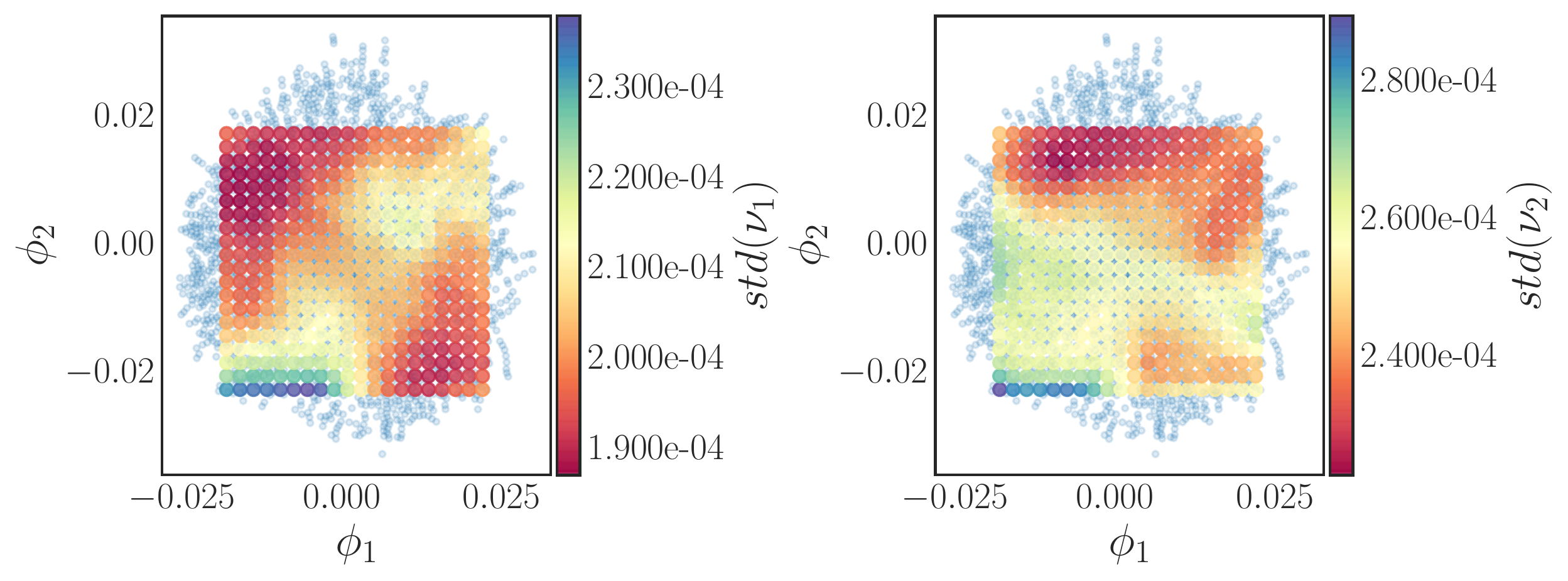}

\includegraphics[scale=0.35]{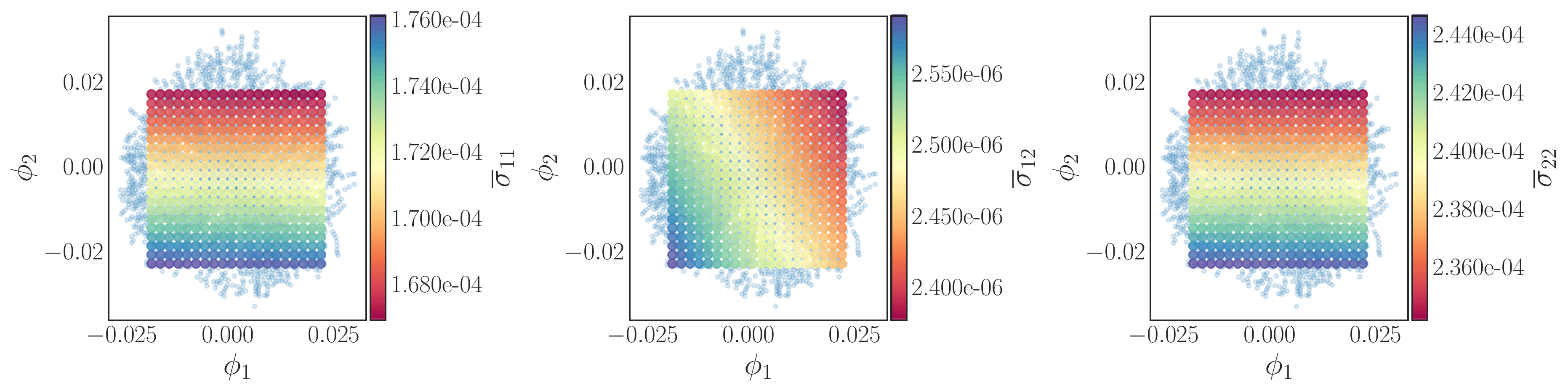}
\includegraphics[scale=0.35]{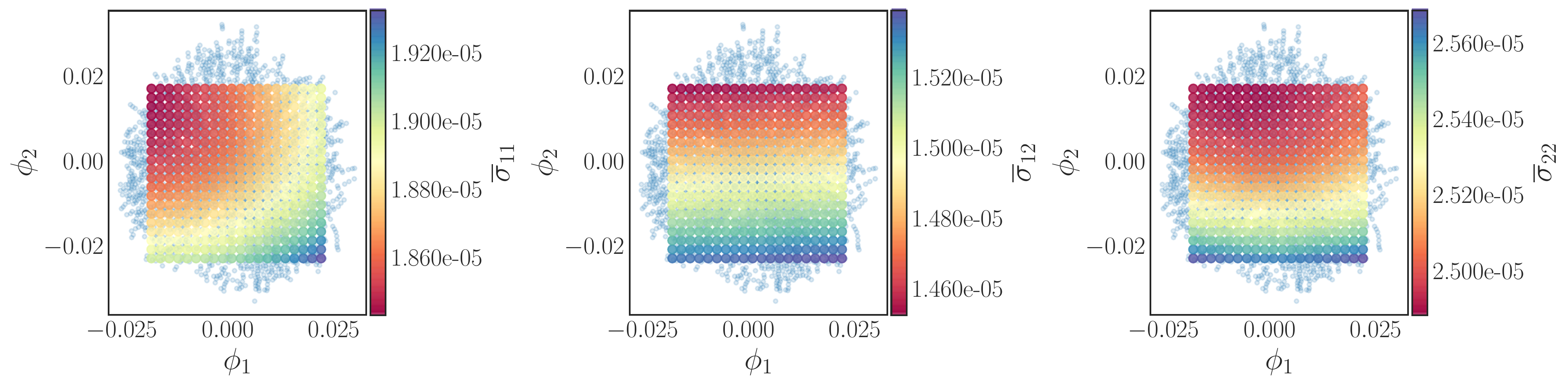}

\caption{[First-Second Rows] illustrate the average drift and the standard deviation of the drift on a test grid (20x20) based on the trained 200 models.
[Third-Fourth Rows] illustrate the average diffusivity and the standard deviation of the diffusivity on a test grid (20x20) based on the trained 200 models. }
\label{fig:Grids}
\end{center}
\end{figure}

In Figure \ref{fig:Grids} we illustrate the same estimations for the 200 models on a 2D grid: we plot the observed values for the mean drift, mean diffusivity based on the 200 models and their standard deviation. The observations that we can draw based on those models are similar to the ones discussed above.

\subsection{Kramers-Moyal vs Neural Network Comparison}

\begin{figure}[ht]
\centering
\begin{center}
\includegraphics[scale=0.35]{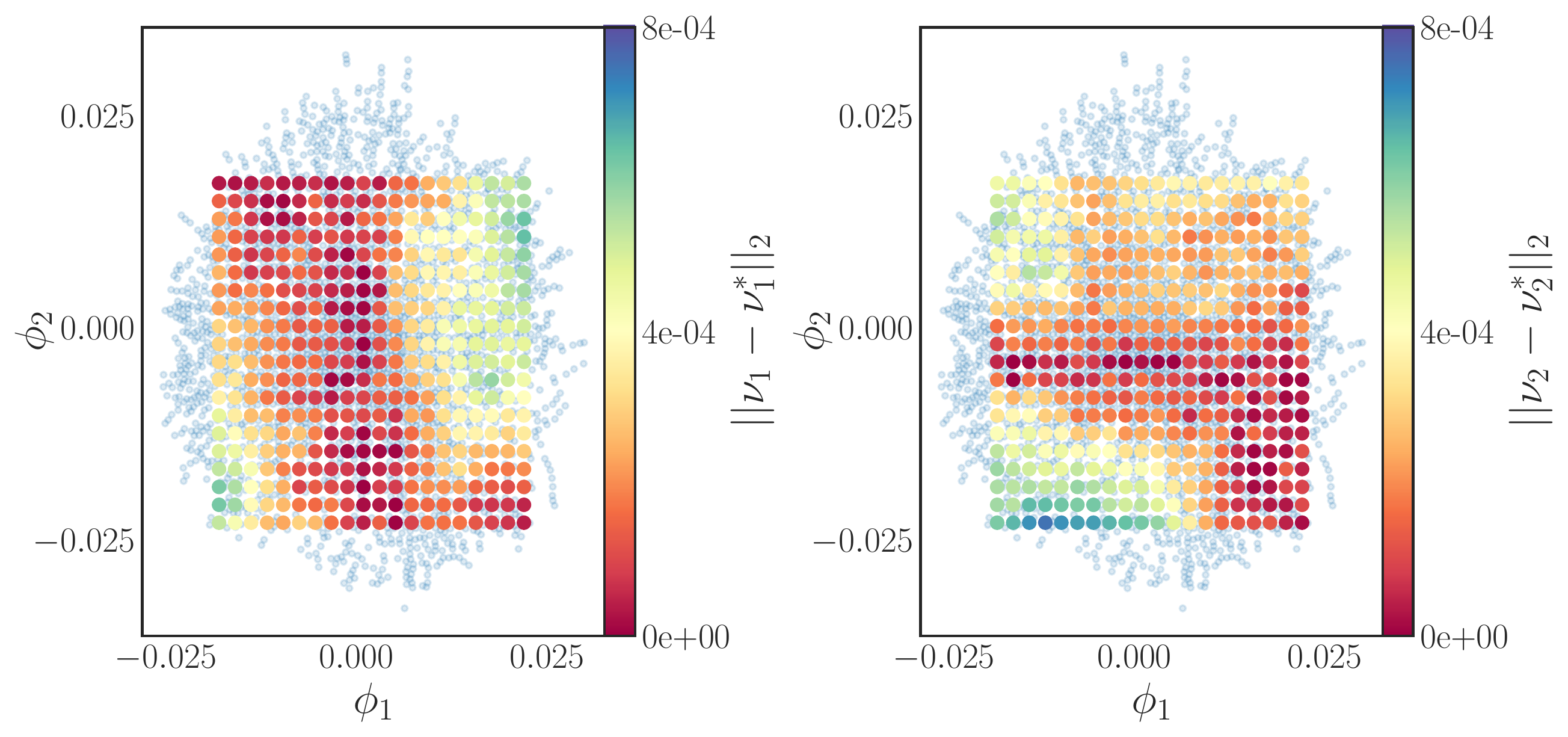}
\includegraphics[scale=0.35 ]{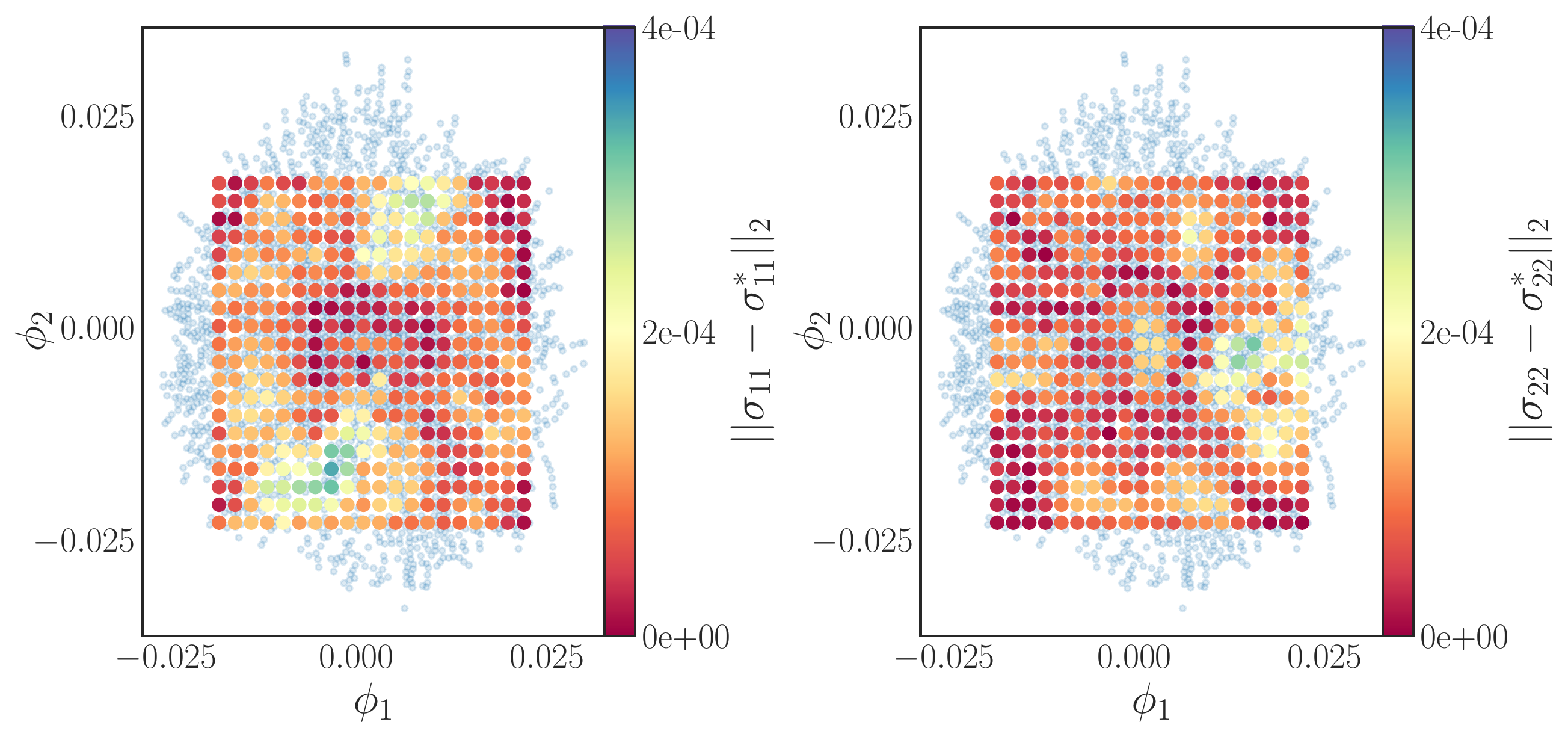}
\caption{[First row] The $l^2$ norm of the difference between the evaluated drift of the network, and by the Kramers-Moyal surrogate model, colored as a function at the nodes of the grid.
[Second Row]The $l^2$ norm of the difference between the evaluated diffusivity of the network and the one from the Kramers-Moyal surrogate model, colored as a function at the nodes of the grid.}
\label{fig:KM_NN_Quantitative_Comparison}
\end{center}
\end{figure}

At the nodes of a Cartesian grid (20x20) we evaluate the drift and diffusivity of the two surrogate models (the Kramers-Moyal and the Neural Network). We then compute the $l^2$ norm between the estimated values of the coefficients with the two models. In Figure \ref{fig:KM_NN_Quantitative_Comparison} we illustrate the grid points (in Diffusion Maps coordinates) colored with the $l^2$ norm for each coefficient. In Table 1 we also provide the mean values of the $l^2$ norms for the drift and the diffusivity.

\begin{table}[th]
 \caption{Estimated mean $l^2$ differences between the estimated drift and diffusivity by the neural network and the Kramers-Moyal surrogate models on a (20 x 20) grid.}
\begin{center}
\begin{tabular}{|c|c|c|c|}
\toprule
Mean $||\nu_1 - \nu_1^* ||_2$    & Mean $||\nu_2 - \nu_2^* ||_2$  & Mean $||\sigma_{11} - \sigma_{11}^* ||_2$ &  Mean $||\sigma_{22} - \sigma_{22}^*||_2$ \\
\midrule
0.000234 & 0.000260 & 0.000104 & 0.000083 \\
\bottomrule
\end{tabular}
\end{center}
\end{table}

The mean values of the $l^2$ norm for the drift and diffusivity reported in the table above are of the same order of magnitude as the standard deviation in Figures \ref{fig:drift_avg_dmaps}, \ref{fig:diffusivity_avg}. This suggests that the discrepancy between the two approaches Kramers-Moyal and neural network falls into the range of uncertainty estimations of the neural network model discussed in Section \ref{sec:Uncertainty_Quantification_NN} of the SI.

\end{document}